\documentclass[12pt,twoside,a4paper]{article}
\usepackage{amssymb}
\usepackage{bbm}
\usepackage{stmaryrd}
\usepackage{amsfonts}
\usepackage{amsmath,amssymb}
\usepackage{mathrsfs}
\usepackage{hyperref}
\usepackage{graphics}
\usepackage{graphicx}
\usepackage{mathrsfs}
\usepackage{epsfig}

\newtheorem{thm}{Theorem}[section]
\newtheorem{cor}[thm]{Corollary}
\newtheorem{lem}[thm]{Lemma}
\newtheorem{prop}[thm]{Proposition}
\newtheorem{defn}[thm]{Definition}

\newtheorem{rem}[thm]{Remark}

\numberwithin{equation}{section}\allowdisplaybreaks

\def\le{\leqslant}
\def\ge{\geqslant}
\def\leq{\leqslant}
\def\geq{\geqslant}

\begin{document}

\begin{center}
\Large\bf Trace Operators for Modulation, $\alpha$-Modulation and
Besov Spaces

\normalsize \vspace*{0.5cm} Hans G.
Feichtinger$^{a,}$\footnote{Email: hans.feichtinger@univie.ac.at},
Chunyan Huang$^{b,}$\footnote{Email:
  chunyanh80@gmail.com},\; Baoxiang Wang$^{b,}$\footnote{Email:
  wbx@math.pku.edu.cn}\\

  \medskip

\footnotesize {\it $^a$Faculty of Mathematics, University of Vienna,
Vienna A-1090, Austria
 .}\\
{\it $^b$LMAM, School of Mathematical Sciences, Peking University,
Beijing 100871, China.}

\vspace*{0.5cm}

\begin{minipage}{13.5cm}
\footnotesize \bf Abstract. \rm  In this paper, we consider the
trace theorem for modulation spaces $M^s_{p,q}$, $\alpha$-modulation
spaces $M^{s,\alpha}_{p,q}$ and Besov spaces $B^s_{p,q}$.  For the
modulation space, we obtain the sharp results.\\

\medskip


\end{minipage}
\end{center}

\section{Introduction} \label{sect1}

\medskip

The $\alpha$-modulation spaces $M^{s,\alpha}_{p,q}$, introduced by
Gr\"obner in \cite{PG} are a class of function spaces that contain
Besov spaces  $B^{s}_{p,q}$ $(\alpha=1)$ and modulation spaces
$M^{s}_{p,q}$ $(\alpha=0)$ as special cases.

There are two kinds of basic coverings on Euclidean $\mathbb {R}^n$
which is very useful in the theory of function spaces and their
applications, one is the uniform covering $\mathbb {R}^n=
\bigcup_{k\in \mathbb{Z}^n} Q_k$, where $Q_k$ denote the unit cube
with center $k$; another is the dyadic covering $\mathbb {R}^n=
\bigcup_{k\in \mathbb{N}} \{\xi: 2^{k-1}\le |\xi| < 2^k\} \bigcup
\{\xi:  |\xi| \le 1\}$. Roughly speaking, these decompositions
together with the frequency-localized techniques yield the
frequency-uniform decomposition operator $\Box_k\sim
\mathscr{F}^{-1} \chi_{Q_k}\mathscr{F}$ and the dyadic decomposition
operator $\Delta_k\sim \mathscr{F}^{-1} \chi_{\{\xi: |\xi|\sim
2^k\}}\mathscr{F}$, respectively. The tempered distributions acted
on these decomposition operators and equipped with the
$\ell^q(L^p(\mathbb{R}^n))$ norms, we then obtain Feichtinger's
modulation spaces and Besov spaces, respectively.

During the past twenty years, the third covering was independently
found by Feichtinger and Gr\"obner \cite{ FeG, Fei1, PG}, and
P\"aiv\"arinta and Somersalo \cite{PS}. This covering, so called
$\alpha$-covering has a moderate scale which is rougher than that of
the uniform covering and  is thinner than that of the dyadic
covering. Applying the $\alpha$-covering to the frequency spaces, in
a similar way as the definition of Besov spaces, Gr\"obner \cite{PG}
introduced the notion of $\alpha$-modulation spaces.

Let $n\ge 2$. For any $x=(x_1,...,x_n) \in \mathbb{R}^n$, we denote
$\bar{x}=(x_1,...,x_{n-1})$. Given a Banach function space
$X(\mathbb{R}^n)$ defined on $\mathbb{R}^n$ and $f\in X$, we ask for
the trace of $f$ on the hyperplane $\{x:\, x=(\bar{x},0)\}$. For the
sake of convenience, this hyperplane will be written as
$\mathbb{R}^{n-1}$. It is clear that a clarification of this problem
is of importance for the boundary value problems of the partial
differential equations. Now we  exactly describe the trace
operators.

\begin{defn}
Let $X$ and $Y$ be quasi-Banach function spaces defined on
$\mathbb{R}^n$ and $\mathbb{R}^{n-1}$, respectively. Denote
\begin{align}
\mathbb{T}\mathbbm{r}: \; f(x)  \to f(\bar{x},0).
\end{align}
If $\mathbb{T}\mathbbm{r}: X\to Y$ and there exists a constant $C>0$
such that
\begin{align}
\|\mathbb{T}\mathbbm{r} f(x)\|_Y \le C  \|f\|_X, \quad \forall f\in
X,
\end{align}
and there exist a continuous linear operator
$\mathbb{T}\mathbbm{r}^{-1}: Y\to X$ such that $
\mathbb{T}\mathbbm{r}\mathbb{T} \mathbbm{r}^{-1}$ is identical
operator, then $\mathbb{T}\mathbbm{r}$ is said to be a retraction
from $X$ onto $Y$.
\end{defn}

If $\mathbb{T}\mathbbm{r}$ is a retraction from $Y$ onto $X$, we see
that the trace of $f\in X$ is well behaved in $Y$. The trace
theorems in modulation spaces and Besov spaces have been extensively
studied. Feichtinger \cite{Fei2} considered the trace theorem for
the modulation space $M^s_{p,q}$ in the case $1\le p,q\le \infty$,
$s>1/q'$ and he obtained that $\mathbb{T}\mathbbm{r}
M^s_{p,q}(\mathbb{R}^n) = M^{s-1/q'}_{p,q} (\mathbb{R}^{n-1})$.
Frazier and Jawerth \cite{FJ} proved that $\mathbb{T}\mathbbm{r}
B^s_{p,q}(\mathbb{R}^n) = B^{s-1/p}_{p,q} (\mathbb{R}^{n-1})$ in the
case $0<p,q\le \infty$ and $s-1/p > \max( (n-1)(1/p-1), 0)$.

\noindent In this paper, we will show the following:

\begin{thm}\label{Th1}
Let $n\geq 2$, $0<p, q\leq \infty$,  $s\in \mathbb{R}$. Then
\begin{align}
\mathbb{T}\mathbbm{r}: f(x) \longrightarrow f(\bar{x}, 0), \;
\bar{x}=(x_1, \cdots, x_{n-1}) \label{tr}
\end{align}
is a retraction from $M_{p, q, \, p\wedge q \wedge
1}^{s}(\mathbb{R}^n)$ onto $M_{p, q}^{s}(\mathbb{R}^{n-1})$.
\end{thm}
Theorem \ref{Th1} is sharp in the sense that $\mathbb{T}\mathbbm{r}:
M_{p, q, \, r}^{s}(\mathbb{R}^n) \not\to  M_{p,
q}^{s}(\mathbb{R}^{n-1})$ for some $r> 1, \; p, \,q \ge 1$. In view
of the basic embedding $M^{s}_{p,q} \subset M^{s_2}_{p,q,q_2}$ for
$s-s_2> 1/q -1/q_2>0$, $s\ge 0$, we immediately have

\begin{cor}\label{cor1}
Let $n\geq 2$, $0<p, q\leq \infty$,  $s\ge 0$. Let
$\mathbb{T}\mathbbm{r}$ be as in \eqref{tr}. Then for any
$\varepsilon>0$,
$$
\mathbb{T}\mathbbm{r}: M_{p, q}^{s+ \frac{1}{p\wedge q\wedge
1}-\frac{1}{q}+ \varepsilon}(\mathbb{R}^n) \to M_{p,
q}^{s}(\mathbb{R}^{n-1}).
$$
\end{cor}
One may ask if Corollary \ref{cor1} holds for the limit case
$s=\varepsilon=0$, we can give a counterexample to show that
$\mathbb{T}\mathbbm{r}: M_{p, q}^{1/q'}(\mathbb{R}^n) \not\to M_{p,
q}^{0}(\mathbb{R}^{n-1})$ in the case $p, \,q>1$.  Write
$$
s_{p} = (n-1)(1/(p\wedge 1) -1).
$$
It is easy to see that $s_p=0$ for $p\ge 1$ and $s_p=(n-1)(1/p-1)$
for $p<1$. For the trace of $\alpha$-modulation spaces, we have the
following result.

\begin{thm}\label{Th2}
Let $n\geq 2$, $0<p, q\leq \infty$,  $s\ge \alpha(n-1)/q + \alpha
s_p$. Let $\mathbb{T}\mathbbm{r}$ be as in \eqref{tr}. Then
$$
\mathbb{T}\mathbbm{r}: M_{p, p\wedge q \wedge 1}^{s+
\alpha/p}(\mathbb{R}^n) \to M_{p, q}^{s}(\mathbb{R}^{n-1}).
$$
\end{thm}

The case $s< \alpha(n-1)/q + \alpha s_p$ is more complicated. We
have the following

\begin{rem}\label{Rem5}
Let $n\geq 2$, $0<p, q\leq \infty$, $s<\alpha (n-1)/q + \alpha s_p$.
Let $\mathbb{T}\mathbbm{r}$ be as in \eqref{tr}. Then
$$
\mathbb{T}\mathbbm{r}: M_{p, p\wedge q \wedge 1}^{s+
\sigma_{\alpha,p,q}}(\mathbb{R}^n) \to M_{p,
q}^{s}(\mathbb{R}^{n-1}),
$$
where
\begin{align}
\sigma_{\alpha, p,q} = \left\{\begin{array}{ll}
\alpha/p + (1-\alpha)[\alpha(n-1)/q +\alpha s_p -s], &  qs +(n-1)(1-\alpha) - q\alpha s_p >0,\\
\alpha/p+ \alpha s_p-s+\varepsilon, &   qs  +(n-1)(1-\alpha) - q\alpha s_p=0,\\
\alpha/p + \alpha s_p -s, &   qs +(n-1)(1-\alpha) - q\alpha s_p <0.
\end{array}\right. \nonumber
\end{align}
\end{rem}

Theorem \ref{Th2} is sharp in the case $s\ge 0$, $p=q=1$. As the end
of this paper, we consider the trace of Besov spaces. If $s>s_p$,
the corresponding result has been obtained in \cite{FJ}. If $s\le
s_p$, we have the following trace theorem for Besov spaces:

\medskip

\begin{thm}\label{Th5}
Let $n\geq 2$, $0<p, q\leq \infty$,  $s\leq s_p$. Let
$\mathbb{T}\mathbbm{r}$ be as in \eqref{tr}. Then we have
$$
\mathbb{T}\mathbbm{r}: \tilde{B}_{ p,p\wedge q \wedge
1}^{s_p+1/p}(\mathbb{R}^n) \to B_{p, q}^{s_p}(\mathbb{R}^{n-1}),
$$
and
$$
\mathbb{T}\mathbbm{r}: B_{ p,p\wedge q \wedge
1}^{s_p+1/p}(\mathbb{R}^n) \to B_{p, q}^{s}(\mathbb{R}^{n-1})
$$
 in the case $s<s_p$. Moreover, when $1<p<\infty$, we have
\begin{align*}
&\mathbb{T}\mathbbm{r}: \tilde{B}_{ p,q\wedge
1}^{1/p,1/p}(\mathbb{R}^n) \to B_{p, q}^{0}(\mathbb{R}^{n-1}),\\
&\mathbb{T}\mathbbm{r}: B_{ p, q \wedge 1}^{1/p,
s+1/p}(\mathbb{R}^n) \to B_{p, q}^{s}(\mathbb{R}^{n-1}), \quad
\text{} \,\,\,s<0.
\end{align*}
\end{thm}
\medskip

The following are some notations which will be frequently used in
this paper: $\mathbb{R}, \mathbb{N}$ and $ \mathbb{Z}$ will stand
for the sets of reals, positive integers and integers, respectively.
$\mathbb{R}_+=[0,\infty)$, $\mathbb{Z}_+= \mathbb{N}\cup \{0\}$.
$c<1$, $C>1$ will denote positive universal constants, which can be
different at different places. $a\lesssim b$ stands for $a\le C b$
for some constant $C>1$, $a\sim b$ means that $a\lesssim b$ and
$b\lesssim a$. We write $a\wedge b =\min(a,b)$, $a\vee b
=\max(a,b)$. We denote by $p'$ the dual number of $p \in
[1,\infty]$, i.e., $1/p+1/p'=1$. We will use Lebesgue spaces
$L^p:=L^p(\mathbb{R}^n)$, $\|\cdot\|_p :=\|\cdot\|_{L^p}$,  We
denote by $\mathscr{S}:=\mathscr{S}(\mathbb{R}^n)$ and
$\mathscr{S}':=\mathscr{S}'(\mathbb{R}^n)$ the Schwartz space and
tempered distribution space, respectively. $B(x,R)$ stands for the
ball in $\mathbb{R}^n$ with center $x$ and radius $R$, $Q(x,R)$
denote the cube in $\mathbb{R}^n$ with center $x$ and side-length
$2R$. $\mathscr{F}$ or \;$\widehat{}$ \;denotes the Fourier
transform; $\mathscr{F}^{-1}$ denotes the inverse Fourier transform.
For any set $A$ with finite elements, we denote by $\# A$ the number
of the elements of $A$.

\section{$\alpha$-modulation spaces}

\subsection{Definition}

A countable set $\mathcal {Q}$ of subsets $Q\subset \mathbb{R}^n$ is
said to be an admissible covering if $\mathbb{R}^n = \bigcup_{Q\in
\mathcal {Q}} Q$ and there exists $n_0 <\infty$ such that $\# \{Q'
\in \mathcal {Q}: \, Q\cap Q'\not= \varnothing\}\le n_0$ for all
$Q\in \mathcal {Q}$. Denote
\begin{align}
& r_Q =\sup \{r\in \mathbb{R}: \, B(c_r, r) \subset Q\}, \nonumber\\
& R_Q =\inf \{R\in \mathbb{R}: \, Q \subset B(c_R,  R)  \}.
\end{align}
Let $0\le \alpha\le 1$. An admissible covering is called an
$\alpha$-covering of $\mathbb{R}^n$, if $|Q| \sim \langle
x\rangle^{\alpha n} $ (uniformly) holds for all $Q\in \mathcal {Q}$
and for all $x\in Q$, and $\sup_{Q\in \mathcal {Q}} R_Q/r_Q \le K$
for some $K<\infty$.

Let $\mathcal {Q}$ be an $\alpha$-covering of $\mathbb{R}^n$. A
corresponding bounded admissible partition of unity of order $p$
($p$-BAPU) $\{\psi_Q\}_{Q\in \mathcal {Q}}$ is a family of smooth
functions satisfying
$$ \left\{
\begin{array}{l}
\psi_Q:  \mathbb {R}^n \to [0,1], \quad {\rm supp} \psi_Q \subset Q,\\
\sum_{Q\in \mathcal {Q}} \psi_Q (\xi) \equiv 1  \quad \forall \xi \in \mathbb{R}^n, \\
\sup_{Q\in \mathcal {Q}} |Q|^{1/(p\wedge 1)-1} \|\mathscr{F}^{-1}
\psi_Q\|_{L^{(p\wedge 1)}} <\infty.
\end{array}
\right.
$$

\begin{defn} \label{alpha-mod-def}
Let $0<p,q\le \infty$, $s\in \mathbb{R}$, $0\le \alpha\le 1$. Let
$\mathcal {Q}$ be an $\alpha$-covering of $\mathbb{R}^n$ with the
$p$-BAPU $\{\psi_Q\}_{Q\in \mathcal {Q}}$. We denote by
$M^{s,\alpha}_{p,q}$ the space of all tempered distributions $f$ for
which the following is finite:
$$
\|f\|_{M^{s, \alpha}_{p,q}}=\left(\sum_{Q\in \mathcal {Q}}\langle
\xi_Q \rangle^{qs} \|\mathscr{F}^{-1}\psi_Q
\mathscr{F}f\|^q_{L^p(\mathbb{R}^n)}\right)^{1/q},
$$
where $\xi_Q\in {Q}$ is arbitrary. For $q=\infty$, we have a usual
substitution for the $\ell^q$ norm with the $\ell^\infty$ norm.
\end{defn}

We now give an exact equivalent norm on $M^{s,\alpha}_{p,q}$. Denote
$$
Q_k=Q (|k|^{\frac{\alpha}{1-\alpha}}k, \, r \langle
k\rangle^{\frac{\alpha}{1-\alpha}}),\quad k\in \mathbb{Z}^n .
$$
It is known that, there exists a constant $r_1>0$ such that for any
$r>r_1$, $\{Q_k\}_{k\in \mathbbm{Z}^n}$ is an $\alpha$-covering of
$\mathbb{R}^n$, i.e., $\mathbb{R}^n=\bigcup_{k\in \mathbb{Z}^n}Q_k$
and there exists $n_0\in \mathbb{N}$ such that $ \#\{l\in
\mathbb{Z}^n: Q_k\cap Q_{k+l}\neq \emptyset  \}\leqslant n_0$.
Moreover, $|Q_k|\sim \langle k\rangle^{\frac{n\alpha}{1-\alpha}}$.
Let $\eta: \mathbb{R} \to [0,1]$ be a smooth bump function
satisfying
\begin{align}
\eta(\xi):= \left\{
\begin{array}{ll}
1, & |\xi|\leq 1,\\
{\rm smooth},  & 1<|\xi| \le 2,\\
 0, & |\xi|\geq 2.
\end{array}
\right.\label{eta}
 \end{align}
We write for $k=(k_1,...,k_n)$ and $\xi=(\xi_1,...,\xi_n)$,
$$
\phi_{k_i}(\xi_i)=\eta\left(\frac{\xi_i-|k|^{\frac{\alpha}{1-\alpha}}k_i}{r\langle
k\rangle^{\frac{\alpha}{1-\alpha}}}\right).
$$
Put
\begin{align} \label{psi-k}
\psi_k(\xi) = \frac{\phi_{k_1}(\xi_1)...
\phi_{k_n}(\xi_n)}{\sum_{k\in \mathbb{Z}^n} \phi_{k_1}(\xi_1)...
\phi_{k_n}(\xi_n) },  \quad k\in \mathbb{Z}^n.
 \end{align}
We have
\begin{lem} \label{lem1.1}
Let $0\leqslant \alpha <1$, $0<p \le \infty$ and $\{\psi_k\}_{k\in
\mathbb{Z}^n}$ be as in \eqref{psi-k}.  Then $\{\psi_k\}_{k\in
\mathbb{Z}^n}$ is a $p$-BAPU for $r>r_1$. In the case $\alpha=0$, we
can take $r_1=1/2$.
\end{lem}

\begin{prop} \label{prop1.1}
Let $0\leqslant \alpha <1$, $0<p, q \leqslant \infty$, then
$$
\|f\|_{M^{s, \alpha}_{p,q}}=\left(\sum_{k\in \mathbb{Z}^n}\langle
k\rangle^{\frac{qs}{1-\alpha}}\|\mathscr{F}^{-1}\psi_k\mathscr{F}f\|^q_{L^p(\mathbb{R}^n)}\right)^{1/q}
$$
is an equivalent norm on $\alpha$-modulation space with the usual
modification for $q=\infty$.
\end{prop}
{\bf Proof.} See \cite{BN}.

\subsection{Equivalent norm via a New p-BAPU}

We now construct a new covering, which is of importance for the
proof of Theorem \ref{Th2}. Let $j \in \mathbb{Z} \backslash \{0\}$.
We divide $[-|j|^{\frac{1}{1-\alpha}}, |j|^{\frac{1}{1-\alpha}}]$
into $2|j|/ \{r_1\}=2 N_j$ intervals with equal length:
$$
[-|j|^{\frac{1}{1-\alpha}}, \,
|j|^{\frac{1}{1-\alpha}}]=[r_{j,-N_j}, r_{j, -N_j+1}]\cup ... \cup
[r_{j, N_j-1}, r_{j, N_j}].
$$
Denote
$$
\mathscr{R}=\{r_{j, s}: j\in \mathbb{N}, s=-N_j, \cdots, N_j\}.
$$
We further write
$$
\mathscr{K}^n_j=\{k= (k_1, \cdots, k_n): k_i \in \mathscr{R}, \,
\max_{1\leq i \leq n} |k_i| =|j|^{\frac{1}{1-\alpha}}\}.
$$
For any $k \in \mathscr{K}^n_j$, we write
$$
Q_{k j}=Q(k, r|j|^{\frac{\alpha}{1-\alpha}}), \quad Q_{k 0}=Q(0, 2).
$$
We will write $\mathscr{K}_j= \mathscr{K}^n_j$ if there is no
confusion.

\begin{prop} \label{prop1.2}
There exists $r_1>0$ such that for any $r>r_1$, $\{Q_{k
j}\}_{k\in\mathscr{K}_j, j\in  \mathbb{Z}_+}$ is an
$\alpha$-covering of $\mathbb{R}^n$.
\end{prop}
{\bf Proof.} Let $j\in  \mathbb{N}\cup \{0\}$. We see that there
exists $r_1>0$ such that for any $r>r_1$,
$\{Q(|j|^{\frac{\alpha}{1-\alpha}}j,
r|j|^{\frac{\alpha}{1-\alpha}})\}_j$ is an $\alpha$-covering of
$\mathbb{R}$. Hence we easily see that
$$
\mathbb{R}\subset \cup_{k\in\mathscr{K}_j, j\in
\mathbb{Z}_+}Q_{kj},\quad  |Q_{kj}|\sim
|j|^{\frac{n\alpha}{1-\alpha}}\sim \langle \xi_{
Q_{kj}}\rangle^{n\alpha}, \forall \,  \xi_{Q_{kj}}\in Q_{k_j},
$$
$$ \#\{Q_{k'j'}: Q_{kj}\cap Q_{k'j'}\neq \emptyset  \}\leqslant
n_0<\infty.
$$
Now, on the basis of the $\alpha$-covering constructed above, we
further construct a $p$-BAPU. Let $j$ be fixed. Denote for $i=1,
\cdots, n$,
$$
\phi_{kj}(\xi_i)=\phi\left(\frac{\xi_i-k_i}{r\langle
j\rangle^{\frac{\alpha}{1-\alpha}}}\right) , \quad k=(k_1, \cdots,
k_n) \in \mathscr{K}_j.
$$
$$
\phi_{kj}(\xi)=\phi_{kj}(\xi_1)...\phi_{kj}(\xi_n).
$$
We put
\begin{align}
\psi_{kj}(\xi)=\frac{\phi_{kj}(\xi)}{\sum_{k\in\mathscr{K}_j, j\in
\mathbb{Z}_+}\phi_{kj}(\xi)}. \label{psikj}
\end{align}

\begin{figure}
\begin{center}
\includegraphics[width=16cm,height=9cm]{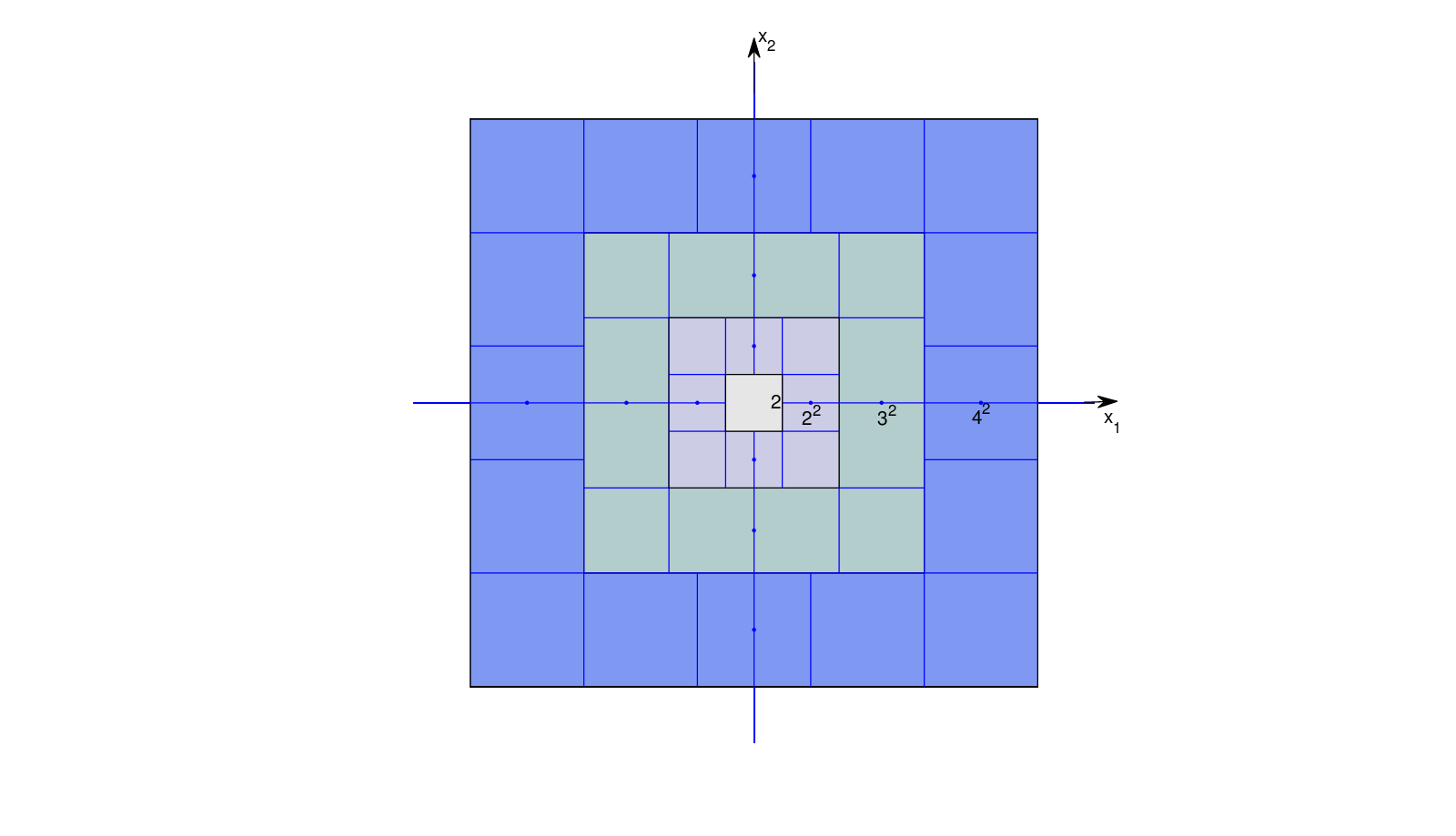}
\begin{minipage}{11cm}
\caption{\footnotesize $\alpha$-covering, the case of $n=2$,
$\alpha=1/2$, $r_1=1$.}
\end{minipage}
\end{center}
\end{figure}

\begin{prop} \label{prop1.2}
Let $0<p <\infty$, $\psi_{kj}$ be as in \eqref{psikj}. Then
$\{\psi_{kj}\}_{k\in \mathscr{K}_j, j\in \mathbb{Z}_+}$ is a
$p$-BAPU.
\end{prop}

Noticing that $|\xi| \sim |j|^{1/(1-\alpha)}$ if $\xi \in
\mathscr{K}_j$, $j\not= 0$,  we immediately have

\begin{prop} \label{prop2.2}
Let $0< \alpha <1$, $0<p, q \leqslant \infty$, then
$$
\|f\|_{M^{s, \alpha}_{p,q}}=\left(\sum_{j\in \mathbb{Z}_+} \langle
j\rangle^{sq/(1-\alpha)} \sum_{k\in \mathscr{K}_j}
\|\mathscr{F}^{-1}\psi_{kj}\mathscr{F}f\|^q_{L^p(\mathbb{R}^n)}\right)^{1/q}
$$
is another equivalent norm on $\alpha$-modulation space.
\end{prop}

\subsection{Modulation spaces}

In the case $\alpha=0$, we get an equivalent norm on modulation
spaces $M^s_{p,q}$:
\begin{align}
\|f\|_{M^{s}_{p,q}}=\left(\sum_{k\in \mathbb{Z}^n}\langle
k\rangle^{qs}\|\mathscr{F}^{-1}\psi_k\mathscr{F}f\|^q_{L^p(\mathbb{R}^n)}\right)^{1/q}.
\label{moddef}
\end{align}
The modulation spaces $M^s_{p,q}$ in the case $0<p,q <1$ was studied
in \cite{Wa1,WH,WHu} by using the norm \eqref{moddef}. Soon after,
Kobayashi \cite{Kub} independently considered such a generalization
in the case $0<p,q<1$.

Recalling that $\bar{x}=(x_1,..., x_{n-1})$, we also define the
following anisotropic modulation spaces $M^s_{p,q,r}$  for which the
norm is defined as
$$
\|f\|_{M^{s}_{p,q,r}}= \left(\sum_{k_n \in \mathbb{Z}}
\left(\sum_{\bar{k}\in \mathbb{Z}^{n-1}}\langle \bar{k}
\rangle^{qs}\|\mathscr{F}^{-1}\psi_k\mathscr{F}f\|^q_{L^p(\mathbb{R}^n)}\right)^{r/q}\right)^{1/r}.
$$
This anisotropic version is of importance for the trace of
modulation spaces.

\subsection{Besov spaces}

Write $\varphi(\cdot)=\eta(\cdot)-\eta(2\cdot)$ and
$\varphi_k:=\varphi(2^{-k}\cdot)$ for $k\geq 1$.
$\varphi_0:=1-\sum_{k\geq 1}\varphi_k$. For simplicity, we write $
\Delta_{k}=\mathscr{F}^{-1}\varphi_{k}\mathscr{F}$. The norm on
Besov spaces $B^s_{p,q}(\mathbb{R}^n)$ are defined as follow:

$$
\|f\|_{B^s_{p,q}}= \left(\sum _{j=0}^{\infty}2^{sjq}\|\Delta_j
f\|^q_p\right)^{1/q}.
$$
For our purpose, we also need the following
\begin{align*}
 \tilde{B}_{ p,q}^{s}(\mathbb{R}^{n})=\left(\sum_{k=0}^{\infty}
 k 2^{skq}\|\Delta_{k} f\|^q_{L^p(\mathbb{R}^{n})}
\right)^{1/q},
\end{align*}

\rm In the case $1<p<\infty$, using Lizorkin's decomposition of
$\mathbb{R}^n$, we have an equivalent quasi-norm on
$B^s_{p,q}(\mathbb{R}^n)$. Let
$$
K_k=\{x: |x_j|<2^k, j=1,2, \ldots, n\} \setminus \{x: |x_j|<2^{k-1},
j=1,2, \ldots, n\}
$$
where $k \in \mathbb{Z}^+$ and
$$
K_0=\{x: |x_j|\leq 1, j=1,2, \ldots, n\}
$$
Subdivide $K_k$ with $k=1,2,3 \ldots,$ by the $3n$ hyper-planes
$\{x:x_m=0\}$ and $\{x: x_m=\pm 2^{k-1}\}$, where $m=1, \ldots, n$,
into cubes $P_{k,t}$. If $k$ is fixed, we obtain $T=4^n-2^n$ cubes.
The cubes near the $n$-th axis are numbered by $t=1,\ldots, 2^n$ in
an arbitrary way and the others are numbered by $t=2^n+1,\ldots, T$.
Let $P_{0,t}=K_0$, if $t=1,\ldots, T$. Then
$$
\mathbb{R}_n=
\cup^{\infty}_{k=0}\bar{K_k}=\cup^{\infty}_{k=0}\cup^{T}_{t=1}\bar{P}_{k,t}.
$$

\begin{figure}
\begin{center}
\includegraphics[width=14cm,height=8cm]{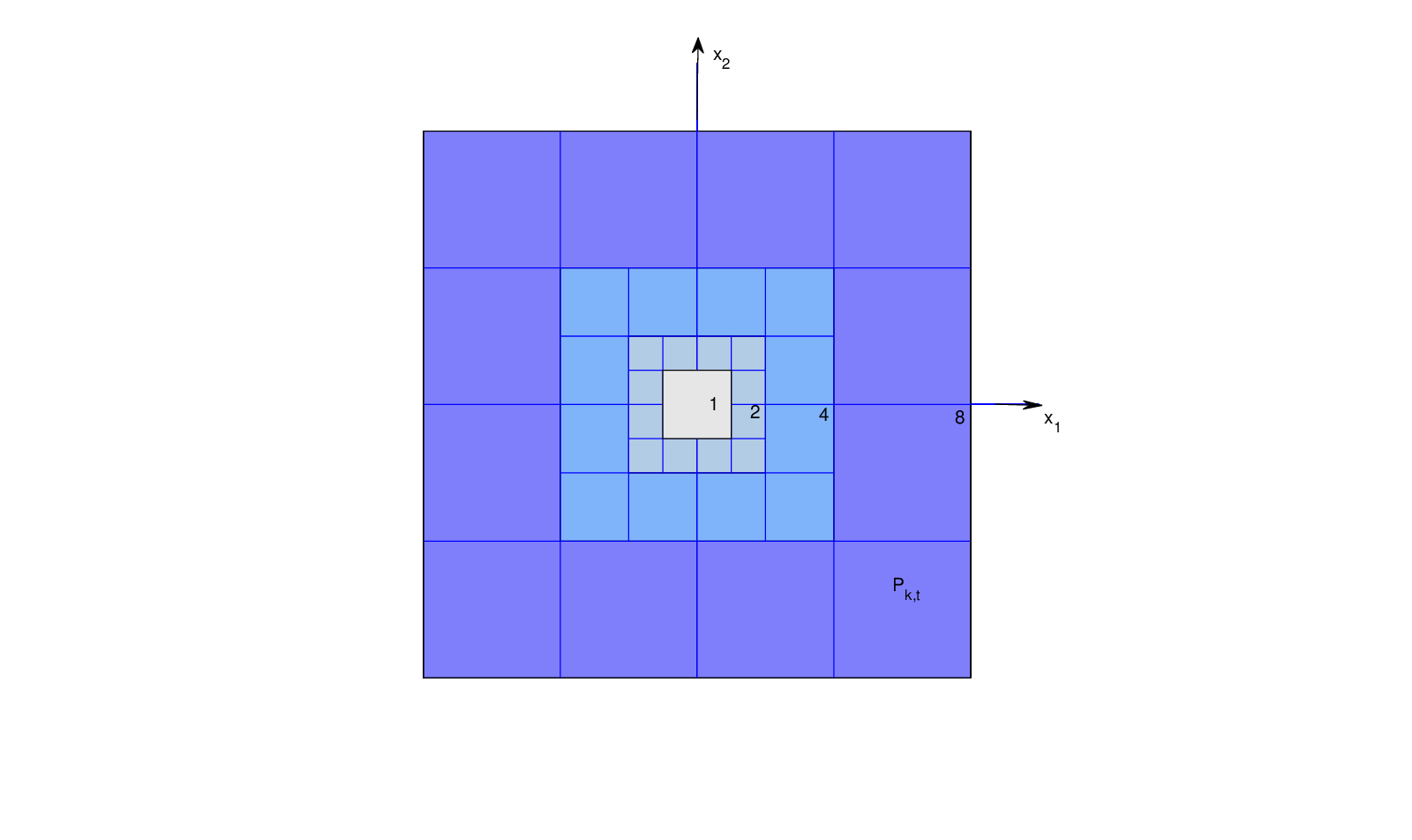}
\begin{minipage}{11cm}
\caption{\footnotesize $1$-covering, the case $n=2$.}
\end{minipage}
\end{center}
\end{figure}

Let $\chi_{k,t}$ be a characteristic function on $P_{k,t}$. Then
$$
\|f\|_{B^s_{p,q}(\mathbb{R}^n)}\asymp \left(\sum_{k=0}^{\infty}
\sum_{t=1}^{T} 2^{skq}\|\mathscr{F}^{-1}\chi_{k,t }\mathscr{F}
f\|^q_{L^p(\mathbb{R}^{n})} \right)^{1/q}.
$$
We construct two new norms. For simplicity, we write
$$
\Delta_{k,t}=\mathscr{F}^{-1}\chi_{k,t}\mathscr{F}.
$$
 Define
\begin{align*}
&B_{ p,q}^{s_1,s_2}=\left(\sum_{k=0}^{\infty}
\left(\sum_{t=1}^{2^{n}} 2^{s_1kq}\|\Delta_{k,t}
f\|^q_{L^p(\mathbb{R}^{n})} + \sum_{t=2^{n}+1}^{T} 2^{s_2
kq}\|\Delta_{k,t}
f\|^q_{L^p(\mathbb{R}^{n})}\right) \right)^{1/q}, \\
& \tilde{B}_{ p,q}^{s_1,s_2}=\left(\sum_{k=0}^{\infty}
\left(\sum_{t=1}^{2^{n}} k 2^{s_1kq}\|\Delta_{k,t}
f\|^q_{L^p(\mathbb{R}^{n})}  + \sum_{t=2^{n}+1}^{T} 2^{s_2
kq}\|\Delta_{k,t} f\|^q_{L^p(\mathbb{R}^{n})}\right)\right)^{1/q}.
\end{align*}

\section{ \bf Proof of Theorem \ref{Th1}}

If there is no explanation, we always assume $r=1/2$ in the $p$-BAPU
for the case of modulation spaces.  To show our main theorem, we
will use the following
\begin{lem}
{\rm (Triebel, \cite{Tr})}\label{Tr1} Let $\Omega$ be a compact
 subset of $\mathbb{R}^n$ and $0<p \leqslant \infty$. Denote
 $L_p^{\Omega}=\{f\in L^p: \textrm{Supp}\hat{f} \subset \Omega\}$. Let
 $0<r<p$. Then
 $$
 \left\|\sup_{z \in \mathbb{R}^n }\frac{f(\cdot - z)}{1+|z|^{n/r}} \right\|_{L^p(\mathbb{R}^n)}
 \lesssim \|f\|_{L^p(\mathbb{R}^n)},
 $$
holds for any $f \in L_p^{\Omega}$.
\end{lem}
Assume that ${\rm Supp} \hat{f}\subset B(\xi_0, R)$. It is easy to
see that for $g=e^{ix \xi_0}f(R^{-1}\cdot)$,
$\hat{g}=R^n\hat{f}(R(\xi-\xi_0))$. It follows that $\textrm{Supp}
\hat{g}\subset B(0, 1)$. Taking $\Omega =B(0, 1)$ in Lemma
\ref{Tr1}, we find that
$$
\left \|\sup_{z \in \mathbb{R}^n }\frac{g(\cdot - z)}{1+|z|^{n/r}}
\right\|_{L^p(\mathbb{R}^n)}
 \lesssim \|g\|_{L^p(\mathbb{R}^n)},
 $$
 By scaling, we have
\begin{equation}
\left \|\sup_{z \in \mathbb{R}^n }\frac{f(\cdot - z)}{1+|Rz|^{n/r}}
\right\|_{L^p(\mathbb{R}^n)}
 \leqslant C\|f\|_{L^p(\mathbb{R}^n)}\label{max}
 \end{equation}
 Note that the constant $C$ in \eqref{max} is independent of $f \in L^p_{B(\xi_0,
 R)}=\{f\in L^p: \textrm{Supp}\hat{f} \subset B(\xi_0,
 R)\}$. It is also independent of $\xi_0 \in \mathbb{R}^n$.

For convenience, we write
$$
\Box_k = \mathscr{F}^{-1}\psi_k \mathscr{F}, \quad k \in
\mathbb{Z}^n.
$$
We define the maximum function $M_k^{\ast} f$ as follows:
\begin{equation}
M_k^{\ast}f= \sup_{y\in \mathbb{Z}^n}\frac{|\Box_k f(x-y)|}{1+|
y|^{n/r}}. \label{max2}
\end{equation}
Taking $y_1=... =y_{n-1}=0$, $y_n=x_n$ in \eqref{max2}, we have for
$ |x_n|\leq 1$,
$$
|(\Box_k f)(\bar{x}, 0)|\lesssim |M_k^{\ast}f(x)|, \quad
\bar{x}=(x_1, \cdots, x_{n-1})
$$
Hence
\begin{equation}
\|(\Box_k f)(\cdot, 0)\|_{L^p(\mathbb{R}^{n-1})}\lesssim
\|M_k^{\ast}f(\cdot, x_n)\|_{L^p(\mathbb{R}^{n-1})}, \label{max3}
\end{equation}
Integrating \eqref{max3} over $x_n \in [0, 1]$, one has that
\begin{equation*}
\|(\Box_k f)(\cdot, 0)\|^p_{L^p(\mathbb{R}^{n-1})}\lesssim
 \int_{\mathbb{R}}\|M_k^{\ast}f(\cdot,
x_n)\|^p_{L^p(\mathbb{R}^{n-1})}dx_n,
\end{equation*}
Hence
\begin{equation}
\|(\Box_k f)(\cdot, 0)\|_{L^p(\mathbb{R}^{n-1})}\lesssim
\|M_k^{\ast}f\|_{L^p(\mathbb{R}^{n})}.  \label{max4}
\end{equation}
We denote by $\mathscr{F}_{\bar{x}}$
$(\mathscr{F}^{-1}_{\bar{\xi}})$ the partial (inverse) Fourier
transform on $\bar{x}$ $(\bar{\xi})$. Write
$\psi_{\bar{k}}(\bar{x})$ as the $p$-BAPU functions in
$\mathbb{R}^{n-1}$ as in \eqref{psi-k}, i.e.,
\begin{align}
\psi_{\bar{k}}(\bar{\xi})=\frac{\phi_{k_1}(\xi_1)...\phi_{k_{n-1}}(\xi_{n-1})
}{\sum_{\bar{k}\in\mathbb{Z}^{(n-1)}}\phi_{k_1}(\xi_1)\cdots\phi_{k_{n-1}}(\xi_{n-1})},
\quad \bar{k}\in \mathbb{Z}^{(n-1)}. \label{psi-k n-1}
\end{align}
Then we have
\begin{align*}
(\mathscr{F}^{-1}_{\bar{\xi}}\psi_{\bar{k}}\mathscr{F}_{\bar{x}})(\bar{x},
0)&=\sum_{l\in
\mathbb{Z}^n}(\mathscr{F}^{-1}_{\bar{\xi}}\psi_{\bar{k}}\mathscr{F}_{\bar{x}}\mathscr{F}^{-1}\psi_{l}\mathscr{F})(\bar{x},
0)\\
&=\sum_{l\in
\mathbb{Z}^n}(\mathscr{F}^{-1}_{\bar{\xi}}\psi_{\bar{k}})\ast(\mathscr{F}^{-1}\psi_{l}\mathscr{F})(\bar{x},
0)
\end{align*}
From the support property of $\psi_l$ as in \eqref{psi-k}, we find
that
$$
\psi_{\bar{k}}\psi_l=0, \quad \text{if } |\bar{l}-\bar{k}|\geq C.
$$
Hence
$$
(\mathscr{F}^{-1}_{\bar{\xi}}\psi_{\bar{k}}\mathscr{F}_{\bar{x}})(\bar{x},
0)=\sum_{l\in \mathbb{Z}^n, \; |\bar{k}-\bar{l}|\leq C
}(\mathscr{F}^{-1}_{\bar{\xi}}\psi_{\bar{k}})\ast
((\mathscr{F}^{-1}\psi_{l}\mathscr{F})(\cdot, 0))
$$

{\it Case 1.} $1\leq p\leq \infty$.  Using Young's inequality,
\eqref{max} and \eqref{max4}, we obtain
\begin{align*}
&\|\mathscr{F}^{-1}_{\bar{\xi}}\psi_{\bar{k}}\mathscr{F}_{\bar{x}}f(\bar{x},
0)\|_{L^p(\mathbb{R}^{n-1})}\\
&\lesssim \sum_{l\in \mathbb{Z}^n, |\bar{k}-\bar{l}|\leq
C}\|\mathscr{F}^{-1}_{\bar{\xi}}\psi_{\bar{k}}\|_{L^1(\mathbb{R}^{n-1})}
\|\mathscr{F}^{-1}\psi_{l}\mathscr{F}f\|_{L^p(\mathbb{R}^{n-1})}\\
&\lesssim \sum_{l\in \mathbb{Z}^n, |\bar{k}-\bar{l}|\leq
C} \|M_l^{\ast}f\|_{L^p(\mathbb{R}^{n})}\\
&\lesssim \sum_{l\in \mathbb{Z}^n, |\bar{k}-\bar{l}|\leq C} \|\Box_l
f\|_{L^p(\mathbb{R}^{n})}.
\end{align*}
Hence,
\begin{align*}
\|f(\bar{x}, 0)\|_{M^{s}_{p,q}(\mathbb{R}^{n-1})} \lesssim
\left(\sum_{\bar{k}\in \mathbb{Z}^{n-1}}\langle \bar{k} \rangle
^{sq}\left(\sum_{l\in \mathbb{Z}^n, |\bar{k}-\bar{l}|\leq C}
\|\Box_l f\|_{L^p(\mathbb{R}^{n})} \right)^q \right)^{1/q}.
\end{align*}
If $0<q \leq 1$, then
\begin{align*}
\|f(\bar{x}, 0)\|_{M^{s}_{p,q}(\mathbb{R}^{n-1})} &\lesssim
\left(\sum_{l\in \mathbb{Z}^{n}}\sum_{\bar{k}\in
\mathbb{Z}^{n-1}}\langle \bar{k} \rangle ^{sq}
\chi_{(|\bar{k}-\bar{l}|\leq
C)}\|\Box_l f\|_{L^p(\mathbb{R}^{n})}^q \right)^{1/q}\\
&\lesssim \left(\sum_{l\in \mathbb{Z}^{n}} \langle \bar{l} \rangle
^{sq} \|\Box_l f\|_{L^p(\mathbb{R}^{n})}^q \right)^{1/q} =
\|f\|_{M_{p, q,q}^{s}}.
\end{align*}
If $1\leq q \leq \infty$, using Minkowski's inequality together with
H\"older's inequality,
\begin{align*}
\|f(\cdot, 0)\|_{M^{s}_{p,q}(\mathbb{R}^{n-1})} &\lesssim
\left(\sum_{\bar{k}, \bar{l}\in \mathbb{Z}^{n-1}}\langle \bar{k}
\rangle ^{sq} \left(\sum_{l_n \in \mathbb{Z}}
 \chi_{(|\bar{k}-\bar{l}|\leq
C)}\|\Box_l f\|_{L^p(\mathbb{R}^{n})}\right)^q\right)^{1/q}\\
 &\lesssim \sum_{l_n \in
\mathbb{Z}}\left(\sum_{\bar{k}, \bar{l}\in \mathbb{Z}^{n-1}} \langle
\bar{k} \rangle ^{sq}\chi_{(|\bar{k}-\bar{l}|\leq
C)}\|\Box_l f\|_{L^p(\mathbb{R}^{n})}^q\right)^{1/q}\\
 &\lesssim \sum_{l_n \in
\mathbb{Z}}\left(\sum_{\bar{l}\in \mathbb{Z}^{n-1}} \langle \bar{l}
\rangle ^{sq} \|\Box_l f\|_{L^p(\mathbb{R}^{n})}^q\right)^{1/q}\\
&=  \|f\|_{M_{p,q,1}^{s}}.
\end{align*}

To begin with the proof for the case  $0<p<1$, we need the following
lemma:

\begin{lem}\label{convo}
Let $0<p\leq 1$. Suppose that $f, g \in L^p_{B(x_0, R)}$, then there
exists a constant $C>0$ which is independent of $x_0 \in
\mathbb{R}^n$ and $R>0$ such that
$$
\|f\ast g\|_p\leq CR^{n(\frac{1}{p}-1)}\|f\|_p\|g\|_p.
$$
\end{lem}
{\bf Proof.} In the case $f, g \in L^p_{B(0, 1)}$, we have
$$
\|f\ast g\|_p \lesssim \|f\|_p\|g\|_p.
$$
Taking $f_{\lambda}=f(\lambda \cdot)$ and $g_{\lambda}=g(\lambda
\cdot)$, we see that
$$f_{R^{-1}},
g_{R^{-1}} \in L^p_{B(0, 1)}, \quad \textrm{if}\,\, f, g \in
L^p_{B(0, R)}.$$ Hence, for any $f, g \in L^p_{B(0, R)}$,

$$
\|f_{R^{-1}}\ast g_{R^{-1}}\|_p \lesssim
\|f_{R^{-1}}\|_p\|g_{R^{-1}}\|_p.
$$
By scaling, we have
$$
\|f\ast g\|_p \lesssim R^{n(\frac{1}{p}-1)}\|f\|_p\|g\|_p.
$$
By a translation $\widehat{e^{i x_0 x}f}=\hat{f}(\xi -x_0)$, we
immediately have the result, as desired. $\hfill \Box$

\medskip

{\it Case 2.} $0< p<1$. By Lemma \ref{convo}, \eqref{max} and
\eqref{max4},
\begin{align*}
&\|\mathscr{F}^{-1}_{\bar{\xi}}\psi_{\bar{k}}\mathscr{F}_{\bar{x}}f(\bar{x},
0)\|^p_{L^p(\mathbb{R}^{n-1})}\\
&\lesssim \sum_{l\in \mathbb{Z}^n, |\bar{k}-\bar{l}|\leq C}
\|\mathscr{F}^{-1}_{\bar{\xi}}\psi_{\bar{k}}\|^p_{L^p(\mathbb{R}^{n-1})}
\|(\mathscr{F}^{-1}\psi_{l}\mathscr{F}f)(\cdot,
0)\|^p_{L^p(\mathbb{R}^{n-1})}\\
&\lesssim \sum_{l\in \mathbb{Z}^n, |\bar{k}-\bar{l}|\leq
C}\|(\mathscr{F}^{-1}\psi_{l}\mathscr{F}f)(\cdot,
0)\|^p_{L^p(\mathbb{R}^{n-1})}\\
&\lesssim \sum_{l\in \mathbb{Z}^n, |\bar{k}-\bar{l}|\leq C} \|\Box_l
f\|^p_{L^p(\mathbb{R}^{n})}.
\end{align*}
It follows that
\begin{align*}
\|f(\bar{x}, 0)\|_{M^{s, \alpha}_{p, q}}  &\lesssim
\left(\sum_{\bar{k}\in \mathbb{Z}^n}\langle \bar{k} \rangle^{sq
}\left(\sum_{l\in \mathbb{Z}^n} \|\Box_l
f\|^p_{L^p(\mathbb{R}^{n})}\chi_{(|\bar{k}-\bar{l}|\leq
C)}\right)^{q/p}\right)^{1/q}.
\end{align*}
If $q \leq p$, one has that
\begin{align*}
\|f(\bar{x}, 0)\|_{M^{s}_{p, q}(\mathbb{R}^{n-1})}  &\lesssim
\left(\sum_{l\in \mathbb{Z}^n}\sum_{k\in \mathbb{Z}^n}\langle
\bar{k} \rangle^{sq} \|\Box^{\alpha}_l
f\|^q_{L^p(\mathbb{R}^{n})}\chi_{(|\bar{k}-\bar{l}|\leq
C)}\right)^{1/q}\\
&\lesssim \|f\|_{M^{s}_{p, q}}.
\end{align*}
If $q \geq p$, using Minkowski's inequality, we have
\begin{align*}
\|f(\bar{x}, 0)\|_{M^{s}_{p, q}(\mathbb{R}^{n-1})} &\lesssim
\left(\sum_{l_n\in \mathbb{Z}}\left(\sum_{\bar{k}, \bar{l}\in
\mathbb{Z}^{n-1}}\langle \bar{k} \rangle^{sq} \|\Box_l
f\|^q_{L^p(\mathbb{R}^{n})}\chi_{(|\bar{k}-\bar{l}|\leq
C)}\right)^{p/q}\right)^{1/p}\\
&\lesssim \left(\sum_{l_n\in \mathbb{Z}}\left(\sum_{ \bar{l}\in
\mathbb{Z}^{n-1}}\langle \bar{l} \rangle^{sq} \|\Box_l
f\|^q_{L^p(\mathbb{R}^{n})}\right)^{p/q}\right)^{1/p}\\
&= \|f\|_{M^{s}_{p,q, p}}.
\end{align*}
In order to show $\mathbb{T}\mathbbm{r}$ is a retraction, we need to
show the existence of $\mathbb{T}\mathbbm{r}^{-1}$. Let $\eta$ be as
in \eqref{eta} satisfying $(\mathscr{F}^{-1}_{\xi_n}\eta) (0)=1$.
For any $f\in M^s_{p,q}(\mathbb{R}^{n-1})$, we define
$$
g(x)= [(\mathscr{F}^{-1}_{\xi_n}\eta) (x_n)] f(\bar{x}):=
(\mathbb{T}\mathbbm{r}^{-1} f) (x).
$$
It is easy to see that $g(\bar{x},0) =f(\bar{x})$ and $\Box_k g=0$
for $|k_n| \ge 3$. Hence,
\begin{align*}
\|g\|_{M^{s}_{p, q, p\wedge q\wedge 1}(\mathbb{R}^{n})} &\lesssim
\left(\sum_{k_n\in \mathbb{Z}}\left(\sum_{\bar{k} \in
\mathbb{Z}^{n-1}}\langle \bar{k} \rangle^{sq} \|\Box_k
g\|^q_{L^p(\mathbb{R}^{n})} \right)^{p\wedge q\wedge 1/q}\right)^{1/p\wedge q\wedge 1}\\
&=  \sum_{|k_n|\le 2}\left(\sum_{\bar{k} \in
\mathbb{Z}^{n-1}}\langle \bar{k} \rangle^{sq} \|\Box_{\bar{k}}
f\|^q_{L^p(\mathbb{R}^{n})}
\|\mathscr{F}^{-1}_{\xi_n}\eta\|^q_{L^p(\mathbb{R})} \right)^{1/q}\\
& \lesssim \|f\|_{M^{s}_{p,q}(\mathbb{R}^{n-1})}.
\end{align*}
It follows that $\mathbb{T}\mathbbm{r}^{-1}:
M^{s}_{p,q}(\mathbb{R}^{n-1}) \to M^{s}_{p,q, p\wedge q\wedge 1
}(\mathbb{R}^{n})$.  $\hfill\Box$\\

\medskip

As the end of this section, we show that Theorem \ref{Th1} and
Corollary \ref{cor1} are sharp conclusions. First, we show that
$\mathbb{T}\mathbbm{r}: M_{p, q, r }^{0}(\mathbb{R}^n) \not\to M_{p,
q}^{0}(\mathbb{R}^{n-1})$ if $r>1$. Let $\eta$ be as in \eqref{eta},
$f = \mathscr{F}^{-1}(\eta(2\xi_1)...\eta(2\xi_n))$. For
$k=(k_1,...,k_n)$, we denote
$$
F(x)= \sum_{|k_n|\le 2^N} \langle k_n\rangle^{-1} e^{{\rm i}k_n x_n}
f (x).
$$
It is easy to see that
$$
\mathscr{F} F (\xi)= \sum_{|k_n|\le 2^N} \langle k_n\rangle^{-1}
\eta(2\xi_1)...\eta(2(\xi_n-k_n)).
$$
Hence, $\Box_k F=0$ if $\max_{i=1,...,n-1} |k_i| >2$ or $|k_n| >
2^N+1$.  In view of the definition
\begin{align*}
\|F\|_{M^{0}_{p, q, r}(\mathbb{R}^{n})} &\lesssim \left(\sum_{|k_n|
\le 2^N+1}\left(\sum_{|k_i| \le 2, \, 1\le i\le n-1} \|\Box_k
F\|^q_{L^p(\mathbb{R}^{n})} \right)^{r/q}\right)^{1/r}\\
&\lesssim \sum_{|k_i| \le 2, \, 1\le i\le n-1} \left(\sum_{|k_n| \le
2^N+1}  \|\Box_k
F\|^r_{L^p(\mathbb{R}^{n})} \right)^{1/r}\\
&\lesssim  \left(\sum_{|k_n| \le 2^N+1} \langle k_n\rangle^{-r}
\right)^{1/r} \lesssim 1.
\end{align*}
On the other hand, we may assume that
$(\mathscr{F}^{-1}_{\xi_n}\eta(2\cdot))(0)=1$.  We have
\begin{align*}
F(\bar{x},0) = \left(\sum_{|k_n|\le 2^N} \langle
k_n\rangle^{-1}\right)
\mathscr{F}^{-1}_{\bar{\xi}}[\eta(2\xi_1)...\eta(2(\xi_{n-1}))].
\end{align*}
So,
\begin{align*}
\|F\|_{M^{0}_{p, q}(\mathbb{R}^{n-1})} & \gtrsim
\left(\sum_{|k_n|\le 2^N} \langle k_n\rangle^{-1}\right)
\left(\sum_{|k_i| \le 2, \, 1\le i\le n-1}
\|\mathscr{F}^{-1}\psi_{\bar{k}} \eta(2\xi_1)...
\eta(2\xi_{n-1})\|^q_{L^p(\mathbb{R}^{n-1})} \right)^{1/q}\\
& \gtrsim  N.
\end{align*}
Let $N\to \infty$, we have $\mathbb{T}\mathbbm{r}: M_{p, q, r
}^{0}(\mathbb{R}^n) \not\to M_{p, q}^{0}(\mathbb{R}^{n-1})$.

Next, we show that $\mathbb{T}\mathbbm{r}: M^{1/q'}_{p, q
}(\mathbb{R}^n) \not\to M_{p, q}^{0}(\mathbb{R}^{n-1})$ as $q>1$.
For $k=(k_1,...,k_n)$, we denote
$$
F(x)= \sum_{|k_n|\le 2^N} \frac{1}{\langle k_n\rangle \ln\langle
k_n\rangle}
 e^{{\rm i}k_n x_n} f (x).
$$
Similarly as in the above, we have
\begin{align*}
\|F\|_{M^{1/q'}_{p, q}(\mathbb{R}^{n})} &\lesssim \left(\sum_{|k_n|
\le 2^N+1} \sum_{|k_i| \le 2, \, 1\le i\le n-1} \langle
k_n\rangle^{q-1}\|\Box_k
F\|^q_{L^p(\mathbb{R}^{n})} \right)^{1/q}\\
&\lesssim  \left(\sum_{|k_n| \le 2^N+1} \frac{1}{\langle k_n\rangle
\ln^q \langle k_n\rangle} \right)^{1/q} \lesssim 1.
\end{align*}
On the other hand,
\begin{align*}
\|F\|_{M^{0}_{p, q}(\mathbb{R}^{n-1})} & \gtrsim
\left(\sum_{|k_n|\le 2^N} \frac{1}{\langle k_n\rangle \ln \langle
k_n\rangle}\right) \to \infty, \quad  N \to \infty.
\end{align*}

\medskip

\section{ \bf Proof of Theorem \ref{Th2} and Remark \ref{Rem5}}

For convenience, we write
$$
\Box^\alpha_{k,j} = \mathscr{F}^{-1}\psi_{k,j}\mathscr{F}, \quad k
\in \mathscr{K}_j, \; j\in \mathbb{Z}_+.
$$
We define the maximum function $M_{k,j}^* f$ as follows:
\begin{equation}
M_{k,j}^*f= \sup_{y\in \mathbb{Z}^n}\frac{|\Box^\alpha_{k,j}f(x-y)|}
{1+|\langle j \rangle^{\alpha/(1-\alpha)} y|^{n/r}}. \label{amax2}
\end{equation}
Taking $y_1=... =y_{n-1}=0$, $y_n=x_n$ in \eqref{amax2}, we have for
$\langle j\rangle^{-\alpha/(1-\alpha)}\leq |x_n|\leq 2\langle j
\rangle^{-\alpha/(1-\alpha)}$,
$$
|(\Box^\alpha_{k,j} f)(\bar{x}, 0)|\lesssim |M_{k,j}^*f(x)|, \quad
\bar{x}=(x_1, \cdots, x_{n-1}).
$$
Hence
\begin{equation}
\|(\Box^\alpha_{k,j} f)(\cdot, 0)\|_{L^p(\mathbb{R}^{n-1})}\lesssim
\|M_{k,j}^*f(\cdot, x_n)\|_{L^p(\mathbb{R}^{n-1})}, \label{amax3}
\end{equation}
Integrating \eqref{amax3} over $x_n \in [\langle j
\rangle^{-\alpha/(1-\alpha)}, 2\langle j
\rangle^{-\alpha/(1-\alpha)}]$, one has that
\begin{equation*}
\|(\Box^\alpha_{k,j} f)(\cdot,
0)\|^p_{L^p(\mathbb{R}^{n-1})}\lesssim \langle j
\rangle^{\alpha/(1-\alpha)}\int_{\mathbb{R}}\|M_{k,j}^*f(\cdot,
x_n)\|^p_{L^p(\mathbb{R}^{n-1})}dx_n,
\end{equation*}
Hence
\begin{equation}
\|(\Box^\alpha_{k,j} f)(\cdot, 0)\|_{L^p(\mathbb{R}^{n-1})} \lesssim
\langle j
\rangle^{\alpha/p(1-\alpha)}\|M_{k,j}^*f\|_{L^p(\mathbb{R}^{n})}.
\label{amax4}
\end{equation}
We denote by $\mathscr{F}_{\bar{x}}(\mathscr{F}^{-1}_{\bar{\xi}})$
the partial (inverse) Fourier transform on
$\bar{x}=(x_1,...,x_{n-1})$ $(\bar{\xi}=(\xi_1,\cdots, \xi_{n-1}))$.
Write $\psi_{m,l}(\bar{x})$ as the $p$-BAPU functions in
$\mathbb{R}^{n-1}$ as in \eqref{psikj}. So, by the definition,
\begin{align}
\|f\|_{M^{s,\alpha}_{p,q}(\mathbb{R}^{n-1})} = \left(\sum_{l\in
\mathbb{Z}_+} \sum_{m\in \mathscr{K}^{n-1}_l} \langle
l\rangle^{\frac{qs}{1-\alpha}}
\|\mathscr{F}^{-1}_{\bar{\xi}}\psi_{m,
l}(\bar{\xi})\mathscr{F}_{\bar{x}} f)(\cdot,
0)\|^q_{L^p(\mathbb{R}^{n-1})} \right)^{1/q}. \label{al-max-1}
\end{align}
In order to have no confusion, we always denote by $\psi_{m,l}$ the
$p$-BAPU function in $\mathbb{R}^{n-1}$ and by $\psi_{k,j}$ the
$p$-BAPU function in $\mathbb{R}^{n}$.  From the support property of
$\psi_{k,j}$, we find that
\begin{align}
(\mathscr{F}^{-1}_{\bar{\xi}}\psi_{m, l}\mathscr{F}_{\bar{x}}
f)(\bar{x}, 0) &=\sum_{j\ge l-C, \, k\in \mathscr{K}^n_j}
(\mathscr{F}^{-1}_{\bar{\xi}}\psi_{m,l}\mathscr{F}_{\bar{x}}\mathscr{F}^{-1}\psi_{k,
j}\mathscr{F} f)(\bar{x}, 0). \label{al-max-2}
\end{align}
For our purpose we further decompose $\mathscr{K}^n_j$. Denote
$$
\mathscr{K}^n_{j,\lambda} = \{k\in \mathscr{K}^n_j: \, \max_{1\le
i\le n-1} |k_i|= \lambda\}, \quad \lambda= r_{j0}, r_{j1},...,
r_{jN_j}, \quad r_{j0}=0.
$$
We easily see that $\sum_{k\in \mathscr{K}^n_j} = \sum_{\lambda=0,
r_{j1},...,r_{jN_j}} \sum_{\mathscr{K}^n_{j,\lambda}}$, $N_j\sim
\langle j \rangle$.  Now we divide our discussion into the following
four cases.

{\it Case 1.}  $1\leq p\leq \infty$ and $0<q\le 1$. By
\eqref{al-max-1} and \eqref{al-max-2},
\begin{align}
\|f(\cdot, 0)\|^q_{M^{s, \alpha}_{p,q}(\mathbb{R}^{n-1})} &\lesssim
\sum_{j\in \mathbb{Z}_+} \sum_{\lambda=0, r_{j1},...,r_{jN_j}}
\sum_{k\in \mathscr{K}^n_{j,\lambda}}  \sum_{l\le j+C} \sum_{ m \in
\mathscr{K}^{n-1}_l}  \langle l \rangle ^{\frac{sq}{1-\alpha}}
\nonumber\\
 & \quad \quad \times  \|
(\mathscr{F}^{-1}_{\bar{\xi}}\psi_{m,l})\ast(\mathscr{F}^{-1}\psi_{k,
j}\mathscr{F} f)(\cdot, 0)\|^q_{L^p(\mathbb{R}^{n-1})} .
\label{al-max-4}
\end{align}
In order to control \eqref{al-max-4} by $\|f\|_{M^{s+\alpha/p,
\alpha}_{p,q}}$, we need to bound the sum $\sum_{l\le j+C} \sum_{ m
\in \mathscr{K}^{n-1}_l}$.  It is easy to see that for fixed $k,j$,
\begin{align}
\# \{m\in \mathscr{K}^{n-1}_l: \, {\rm supp\, }\psi_{m,l} \cap {\rm
supp\, }\psi_{k,j} (\cdot,0) \not= \varnothing\} \lesssim \min
\left( \langle l\rangle^{n-2},  \, \langle
j\rangle^{\frac{\alpha(n-2)}{1-\alpha}}/ \langle
l\rangle^{\frac{\alpha(n-2)}{1-\alpha}} \right). \label{number}
\end{align}
Moreover, $k\in \mathscr{K}^n_{j, r_{ja}}$ means that ${\rm supp\,
}\psi_{m,l} \cap {\rm supp\, }\psi_{k,j} (\cdot,0) \not=
\varnothing$ only if $a^{1-\alpha} \langle j\rangle^\alpha \lesssim
l \lesssim (1+a)^{1-\alpha} \langle j\rangle^\alpha$. Hence,  in
view of Young's inequality, \eqref{amax4},
\begin{align}
\Delta_{ja} & := \sum_{ k \in \mathscr{K}^n_{j, r_{ja}}} \sum_{m \in
\mathscr{K}^{n-1}_l, \, l\le j+C } \langle l \rangle
^{\frac{sq}{1-\alpha}} \| (\mathscr{F}^{-1}_{\bar{\xi}}\psi_{m,l})*
(\mathscr{F}^{-1}\psi_{k, j}\mathscr{F}) f
(\cdot, 0)\|^q_{L^p(\mathbb{R}^{n-1})} \nonumber\\
& \lesssim  \sum_{ k \in \mathscr{K}^n_{j, r_{ja}}}
\sum_{a^{1-\alpha} \langle j\rangle^\alpha \lesssim l \lesssim
(1+a)^{1-\alpha} \langle j\rangle^\alpha }  \langle l \rangle
^{\frac{sq}{1-\alpha}} \min \left( \langle l\rangle^{n-2},  \,
\langle j\rangle^{\frac{\alpha(n-2)}{1-\alpha}}/ \langle
l\rangle^{\frac{\alpha(n-2)}{1-\alpha}} \right) \nonumber\\
& \quad \quad \times \| (\mathscr{F}^{-1}\psi_{k, j}\mathscr{F}) f
(\cdot, 0)\|^q_{L^p(\mathbb{R}^{n-1})}. \nonumber\\
& \lesssim  \sum_{ k \in \mathscr{K}^n_{j, r_{ja}}}
\sum_{a^{1-\alpha} \langle j\rangle^\alpha \lesssim l \lesssim
(1+a)^{1-\alpha} \langle j\rangle^\alpha }  \langle l \rangle
^{\frac{sq}{1-\alpha}} \min \left( \langle l\rangle^{n-2},  \,
\langle j\rangle^{\frac{\alpha(n-2)}{1-\alpha}}/ \langle
l\rangle^{\frac{\alpha(n-2)}{1-\alpha}} \right)  \nonumber\\
& \quad \quad \times \langle
j\rangle^{\frac{q\alpha}{p(1-\alpha)}}\| \Box^\alpha_{k,j} f
\|^q_{L^p(\mathbb{R}^{n})}. \label{al-max-5}
\end{align}
We discuss the following four subcases.

{\it Case 1A.} $\alpha (n-1) \le qs$. If $a=0$, one has that
\begin{align}
\Delta_{ja}  & \lesssim  \sum_{ k \in \mathscr{K}^n_{j, r_{ja}}}
\sum_{0 \le  l \lesssim  \langle j\rangle^\alpha } \langle l \rangle
^{\frac{sq}{1-\alpha} +(n-2)}  \langle
j\rangle^{\frac{q\alpha}{p(1-\alpha)}}\| \Box^\alpha_{k,j} f
\|^q_{L^p(\mathbb{R}^{n})} \nonumber\\
& \lesssim  \sum_{ k \in \mathscr{K}^n_{j, r_{ja}}} \langle
j\rangle^ {\frac{q}{1-\alpha}\left(\alpha s+
\frac{\alpha}{p}+(1-\alpha) \frac{\alpha(n-1)}{q}\right)}  \|
\Box^\alpha_{k,j} f \|^q_{L^p(\mathbb{R}^{n})}\nonumber\\
& \lesssim  \sum_{ k \in \mathscr{K}^n_{j, r_{ja}}} \langle
j\rangle^ {\frac{q}{1-\alpha}\left(s+ \frac{\alpha}{p}\right)} \|
\Box^\alpha_{k,j} f \|^q_{L^p(\mathbb{R}^{n})}. \label{al-max-5a}
\end{align}
If $a\ge 1$, we have
\begin{align}
\Delta_{ja}  & \lesssim  \sum_{ k \in \mathscr{K}^n_{j, r_{ja}}}
\sum_{a^{1-\alpha} \langle j\rangle^\alpha \lesssim  l \lesssim
(1+a)^{1-\alpha} \langle j\rangle^\alpha } \langle l \rangle
^{\frac{sq}{1-\alpha} - \frac{\alpha (n-2)}{1-\alpha}} \langle
j\rangle^{\frac{q\alpha}{p(1-\alpha)} + \frac{\alpha
(n-2)}{1-\alpha} }\| \Box^\alpha_{k,j} f
\|^q_{L^p(\mathbb{R}^{n})} \nonumber\\
& \lesssim  \sum_{ k \in \mathscr{K}^n_{j, r_{ja}}}
a^{qs-\alpha(n-1)}
  \langle j\rangle^{\frac{\alpha
sq}{1-\alpha}+ \frac{q\alpha}{p(1-\alpha)} + \alpha (n-1) }\|
\Box^\alpha_{k,j} f \|^q_{L^p(\mathbb{R}^{n})}. \label{al-max-5aa}
\end{align}
It follows from $qs\ge \alpha(n-1)$ that $\Delta_{ja}$ takes the
maximal value as $a=N_j \sim \langle j\rangle$. Hence,
\begin{align}
\Delta_{ja} & \lesssim  \sum_{ k \in \mathscr{K}^n_{j, r_{ja}}}
\langle j\rangle^ {\frac{q}{1-\alpha}\left(s+
\frac{\alpha}{p}\right)} \| \Box^\alpha_{k,j} f
\|^q_{L^p(\mathbb{R}^{n})}. \label{al-max-6}
\end{align}
Inserting the estimates of $\Delta_{ja}$ as in \eqref{al-max-5a} and
\eqref{al-max-6} into \eqref{al-max-4}, we have
\begin{align}
\|f(\cdot, 0)\|^q_{M^{s, \alpha}_{p,q}(\mathbb{R}^{n-1})} &\lesssim
\sum_{j\in \mathbb{Z}_+} \sum_{k\in \mathscr{K}^n_{j}}   \langle j
\rangle ^{\frac{q}{1-\alpha}(s+\frac{\alpha}{p})}
\|\Box^\alpha_{k,j} f\|^q_{L^p(\mathbb{R}^{n})} =
\|f\|^q_{M^{s+\alpha/p, \alpha}_{p,q}}. \nonumber
\end{align}
{\it Case 1B.} $\alpha (n-1) > qs$ and $qs+ (1-\alpha) (n-1)>0$. If
$a\ge 1$, from \eqref{al-max-5aa} and $\alpha (n-1) > qs$ we have
\begin{align}
\Delta_{ja} & \lesssim  \sum_{ k \in \mathscr{K}^n_{j, r_{ja}}}
\langle j\rangle^{\frac{\alpha sq}{1-\alpha}+
\frac{q\alpha}{p(1-\alpha)} + \alpha (n-1) }\| \Box^\alpha_{k,j} f
\|^q_{L^p(\mathbb{R}^{n})}. \label{al-max-7}
\end{align}
Since $qs+ (1-\alpha) (n-1)>0$,  similar to \eqref{al-max-5a}, we
see that \eqref{al-max-7} also holds for the case $a=0$. It follows
that
\begin{align}
\|f(\cdot, 0)\|_{M^{s, \alpha}_{p,q}(\mathbb{R}^{n-1})} &\lesssim
 \|f\|_{M^{\alpha s+ (1-\alpha)\frac{\alpha(n-1)}{q}+\alpha/p, \, \alpha}_{p,q}(\mathbb{R}^{n})}. \nonumber
\end{align}

{\it Case 1C.} $qs=-(n-1)(1-\alpha)$. Using the first estimate as in
\eqref{al-max-5a}, we have for $a=0$,
\begin{align}
\Delta_{ja} & \lesssim  \sum_{ k \in \mathscr{K}^n_{j, r_{ja}}}
\langle j\rangle^{ \frac{q\alpha}{p(1-\alpha)} } \ln \langle
j\rangle \| \Box^\alpha_{k,j} f \|^q_{L^p(\mathbb{R}^{n})}.
\label{al-max-9}
\end{align}
For $a\ge 1$,
\begin{align}
\Delta_{ja} & \lesssim  \sum_{ k \in \mathscr{K}^n_{j, r_{ja}}}
\langle j\rangle^{\frac{\alpha sq}{1-\alpha}+
\frac{q\alpha}{p(1-\alpha)} + \alpha (n-1) }\| \Box^\alpha_{k,j} f
\|^q_{L^p(\mathbb{R}^{n})}. \label{al-max-10}
\end{align}
This implies that
\begin{align}
\|f(\cdot, 0)\|_{M^{s, \alpha}_{p,q}(\mathbb{R}^{n-1})} &\lesssim
 \|f\|_{M^{\alpha/p+ \varepsilon, \, \alpha}_{p,q}(\mathbb{R}^{n})}. \nonumber
\end{align}

{\it Case 1D.} $qs<-(n-1)(1-\alpha)$. It is easy to see that $\ln
\langle j\rangle$ can be removed in \eqref{al-max-9}. So, we have
the result, as desired.

\medskip

{\it Case 2.} $1\le p,q\le \infty$. Using Minkowski's inequality, we
have
\begin{align}
& \|f(\cdot, 0)\|_{M^{s, \alpha}_{p,q}(\mathbb{R}^{n-1})}
\nonumber\\
& \lesssim  \left( \sum_{l \in \mathbb{Z}_+} \sum_{ m \in
\mathscr{K}^{n-1}_l} \langle l \rangle ^{\frac{sq}{1-\alpha}} \left(
\sum_{j\in \mathbb{Z}_+}  \sum_{k\in \mathscr{K}^n_{j}} \|
(\mathscr{F}^{-1}_{\bar{\xi}}\psi_{m,l})\ast(\mathscr{F}^{-1}\psi_{k,
j}\mathscr{F} f)(\cdot, 0)\|_{L^p(\mathbb{R}^{n-1})} \right)^q
\right)^{1/q} \nonumber\\
& \lesssim   \sum_{j\in \mathbb{Z}_+} \sum^{N_j}_{a=0} \sum_{k\in
\mathscr{K}^n_{j,r_{ja}}} \left(\sum_{l \in \mathbb{Z}_+} \sum_{ m
\in \mathscr{K}^{n-1}_l} \langle l \rangle ^{\frac{sq}{1-\alpha}}
 \|
(\mathscr{F}^{-1}_{\bar{\xi}}\psi_{m,l})\ast(\mathscr{F}^{-1}\psi_{k,
j}\mathscr{F} f)(\cdot, 0)\|^q_{L^p(\mathbb{R}^{n-1})} \right)^{1/q}
. \label{al-max-14}
\end{align}
Using the same way as in Case 1, we can get that for any $k\in
\mathscr{K}^n_{j,r_{ja}}$,
\begin{align}
&  \sum_{l \in \mathbb{Z}_+} \sum_{ m \in \mathscr{K}^{n-1}_l}
\langle l \rangle ^{\frac{sq}{1-\alpha}}
 \|(\mathscr{F}^{-1}_{\bar{\xi}}\psi_{m,l})\ast(\mathscr{F}^{-1}\psi_{k,
j}\mathscr{F} f)(\cdot, 0)\|^q_{L^p(\mathbb{R}^{n-1})}
\nonumber\\
&  \lesssim \sum_{a^{1-\alpha} \langle j\rangle^\alpha \lesssim l
\lesssim (1+a)^{1-\alpha} \langle j \rangle^\alpha } \langle l
\rangle ^{\frac{sq}{1-\alpha}} \min \left( \langle l\rangle^{n-2},
\, \frac{\langle j\rangle^{\frac{\alpha(n-2)}{1-\alpha}}}{\langle
l\rangle^{\frac{\alpha(n-2)}{1-\alpha}} } \right) \langle
j\rangle^{\frac{q\alpha}{p(1-\alpha)}}\| \Box^\alpha_{k,j} f
\|^q_{L^p(\mathbb{R}^{n})} . \label{al-max-15}
\end{align}
Repeating the calculation procedure as in Case 1, we have the
result, as desired.

\medskip

{\it Case 3}. $0<q\le p<1$. We have
\begin{align}
\|\mathscr{F}^{-1}_{\bar{\xi}}\psi_{m, l}\mathscr{F}_{\bar{x}} f
(\cdot, 0)\|^p_{L^p(\mathbb{R}^{n-1})} \le \sum_{j\ge l-C, \, k\in
\mathscr{K}^n_j}
\|(\mathscr{F}^{-1}_{\bar{\xi}}\psi_{m,l}\mathscr{F}_{\bar{x}}\mathscr{F}^{-1}\psi_{k,
j}\mathscr{F} f)(\cdot, 0)\|^p_{L^p(\mathbb{R}^{n-1})}.
\label{al-max-16}
\end{align}
It follows that
\begin{align}
& \|f(\cdot, 0)\|_{M^{s, \alpha}_{p,q}(\mathbb{R}^{n-1})}
\nonumber\\
& \lesssim  \left( \sum_{l \in \mathbb{Z}_+} \sum_{ m \in
\mathscr{K}^{n-1}_l} \langle l \rangle ^{\frac{sq}{1-\alpha}} \left(
\sum_{j\in \mathbb{Z}_+}  \sum_{k\in \mathscr{K}^n_{j}} \|
(\mathscr{F}^{-1}_{\bar{\xi}}\psi_{m,l})\ast(\mathscr{F}^{-1}\psi_{k,
j}\mathscr{F} f)(\cdot, 0)\|^p_{L^p(\mathbb{R}^{n-1})} \right)^{q/p}
\right)^{1/q} \nonumber\\
& \lesssim  \left( \sum_{j\in \mathbb{Z}_+} \sum^{N_j}_{a=0}
\sum_{k\in \mathscr{K}^n_{j,r_{ja}}} \sum_{l \le j+C} \sum_{ m \in
\mathscr{K}^{n-1}_l} \langle l \rangle ^{\frac{sq}{1-\alpha}}
 \|
(\mathscr{F}^{-1}_{\bar{\xi}}\psi_{m,l})\ast(\mathscr{F}^{-1}\psi_{k,
j}\mathscr{F} f)(\cdot, 0)\|^q_{L^p(\mathbb{R}^{n-1})} \right)^{1/q}
. \label{al-max-17}
\end{align}
For convenience, we write
\begin{align}
\Upsilon_{ja} :=  \sum_{k\in \mathscr{K}^n_{j,r_{ja}}} \sum_{l \le
j+C} \sum_{ m \in \mathscr{K}^{n-1}_l} \langle l \rangle
^{\frac{sq}{1-\alpha}}
 \|
(\mathscr{F}^{-1}_{\bar{\xi}}\psi_{m,l})\ast(\mathscr{F}^{-1}\psi_{k,
j}\mathscr{F} f)(\cdot, 0)\|^q_{L^p(\mathbb{R}^{n-1})} .
\label{al-max-18}
\end{align}
It follows from Lemma \ref{convo}, the property of $p$-BAPU and
\eqref{number} that
\begin{align}
\Upsilon_{ja} & \lesssim  \sum_{k\in \mathscr{K}^n_{j,r_{ja}}}
\sum_{a^{1-\alpha}\langle j\rangle^\alpha \lesssim l \lesssim
(1+a)^{1-\alpha}\langle j\rangle^\alpha} \sum_{ m \in
\mathscr{K}^{n-1}_l} \langle l \rangle ^{\frac{sq}{1-\alpha}}
\nonumber\\
& \quad\quad \times \langle
j\rangle^{q\frac{\alpha(n-1)}{1-\alpha}(\frac{1}{p}-1)} \|
(\mathscr{F}^{-1}_{\bar{\xi}}\psi_{m,l})\|^q_{L^p(\mathbb{R}^{n-1})}
\|(\mathscr{F}^{-1}\psi_{k, j}\mathscr{F} f)(\cdot,
0)\|^q_{L^p(\mathbb{R}^{n-1})} \nonumber\\
 & \lesssim  \sum_{k\in \mathscr{K}^n_{j,r_{ja}}}
\sum_{a^{1-\alpha}\langle j\rangle^\alpha \lesssim l \lesssim
(1+a)^{1-\alpha}\langle j\rangle^\alpha}  \langle l \rangle
^{\frac{sq}{1-\alpha}} \langle
j\rangle^{q\frac{\alpha(n-1)}{1-\alpha}(\frac{1}{p}-1)} \langle
l\rangle^{-q\frac{\alpha(n-1)}{1-\alpha}(\frac{1}{p}-1)}
\nonumber\\
& \quad\quad \times  \min \left( \langle l\rangle^{n-2}, \, \langle
j\rangle^{\frac{\alpha(n-2)}{1-\alpha}} \langle
l\rangle^{-\frac{\alpha(n-2)}{1-\alpha}}  \right)  \langle
j\rangle^{\frac{q\alpha}{p(1-\alpha)}} \|\Box^\alpha_{k, j}
f\|^q_{L^p(\mathbb{R}^{n})} . \label{al-max-20}
\end{align}
If $a\ge 1$, then we have
\begin{align}
\Upsilon_{ja}
 & \lesssim  \sum_{k\in \mathscr{K}^n_{j,r_{ja}}}
a ^{sq -\alpha(n-1)- q\alpha(n-1)(\frac{1}{p}-1)}  \langle
j\rangle^{ \frac{\alpha q s}{1-\alpha} + q
\alpha(n-1)(\frac{1}{p}-1) + \alpha(n-1) +
\frac{q\alpha}{p(1-\alpha)}}
 \|\Box^\alpha_{k, j} f\|^q_{L^p(\mathbb{R}^{n})} .
\label{al-max-21}
\end{align}
If $a=0$,
\begin{align}
\Upsilon_{ja}
 & \lesssim  \sum_{k\in \mathscr{K}^n_{j,r_{ja}}}
\sum_{0 \le l \lesssim \langle j\rangle^\alpha}  \langle l \rangle
^{\frac{sq}{1-\alpha} -q\frac{\alpha(n-1)}{1-\alpha}(\frac{1}{p}-1)
+(n-2) } \langle
j\rangle^{q\frac{\alpha(n-1)}{1-\alpha}(\frac{1}{p}-1) +
\frac{q\alpha}{p(1-\alpha)}}
 \|\Box^\alpha_{k, j} f\|^q_{L^p(\mathbb{R}^{n})} .
\label{al-max-22}
\end{align}
Now we divide our discussion into the following four subcases.

{\it Case 3A}. $qs- \alpha(n-1) - q\alpha(n-1)(1/p-1) \ge 0$. If
$a\ge 1$, we see that the upper bound in \eqref{al-max-21} will be
attained at $a\sim \langle j\rangle$. If $a=0$, the summation on $l$
in \eqref{al-max-22} can be easily controlled. Anyway, we have
\begin{align}
\Upsilon_{ja}
 & \lesssim  \sum_{k\in \mathscr{K}^n_{j,r_{ja}}}
 \langle
j\rangle^{ \frac{\alpha q s}{1-\alpha} +
\frac{q\alpha}{p(1-\alpha)}}
 \|\Box^\alpha_{k, j} f\|^q_{L^p(\mathbb{R}^{n})} .
\label{al-max-23}
\end{align}
Combining \eqref{al-max-17} and \eqref{al-max-23}, we immediately
have the result, as desired.

{\it Case 3B}. $q\alpha s_p-(1-\alpha)(n-1) < qs < q \alpha s_p
+\alpha(n-1)$. In this case, using \eqref{al-max-21} and
\eqref{al-max-22},  we can repeat the procedure as in Case 1B to get
the result and we omit the details of the proof.

{\it Case 3C}. $q\alpha s_p-(1-\alpha)(n-1) =qs $. This case is
similar to Case 1C.

{\it Case 3D}. $q\alpha s_p-(1-\alpha)(n-1) >qs $. We can deal with
this case by following the same way as in Case 1D.

{\it Case 4}. $0< p<1$, $q>p$. By Minkowski's inequality, we have
\begin{align}
& \|f(\cdot, 0)\|_{M^{s, \alpha}_{p,q}(\mathbb{R}^{n-1})}
\nonumber\\
& \lesssim  \left( \sum_{l \in \mathbb{Z}_+} \sum_{ m \in
\mathscr{K}^{n-1}_l} \langle l \rangle ^{\frac{sq}{1-\alpha}} \left(
\sum_{j\in \mathbb{Z}_+}  \sum_{k\in \mathscr{K}^n_{j}} \|
(\mathscr{F}^{-1}_{\bar{\xi}}\psi_{m,l})\ast(\mathscr{F}^{-1}\psi_{k,
j}\mathscr{F} f)(\cdot, 0)\|^p_{L^p(\mathbb{R}^{n-1})} \right)^{q/p}
\right)^{1/q} \nonumber\\
& \lesssim  \left(  \sum_{j\in \mathbb{Z}_+}  \sum_{k\in
\mathscr{K}^n_{j}} \left( \sum_{l \in \mathbb{Z}_+} \sum_{ m \in
\mathscr{K}^{n-1}_l} \langle l \rangle ^{\frac{sq}{1-\alpha}} \|
(\mathscr{F}^{-1}_{\bar{\xi}}\psi_{m,l})\ast(\mathscr{F}^{-1}\psi_{k,
j}\mathscr{F} f)(\cdot, 0)\|^q_{L^p(\mathbb{R}^{n-1})} \right)^{p/q}
\right)^{1/p}. \label{al-max-24}
\end{align}
Then we can repeat the procedures as in the proof of Theorem
\ref{Th1} and the above techniques in Case 3 to have the result, as
desired. The details of the proof are omitted. $\hfill\Box$

\medskip
\section{Proof of Theorem \ref{Th5}}

Now we prove Theorem \ref{Th5}. Now we define the maximum function
$M_{k}^{\ast} f$ as follows:
\begin{equation}
M_{k}^{\ast}f= \sup_{y\in \mathbb{Z}^n}\frac{|\Delta_{k}
f(x-y)|}{1+| 2^{k}y|^{n/r}}.\label{max9}
\end{equation}
Taking $y_1=... =y_{n-1}=0$, $y_n=x_n$ in \eqref{max9}, we have for
$ 2^{-k-1}\leq |x_n|\leq 2^{-k}$,
$$
|(\Delta_{k} f)(\bar{x}, 0)|\lesssim |M_{k}^{\ast}f(x)|, \quad
\bar{x}=(x_1, \cdots, x_{n-1})
$$
Hence
\begin{equation}
\|(\Delta_{k} f)(\cdot, 0)\|_{L^p(\mathbb{R}^{n-1})}\lesssim
\|M_{k}^{\ast}f(\cdot, x_n)\|_{L^p(\mathbb{R}^{n-1})}, \label{max10}
\end{equation}
Integrating \eqref{max10}, one has that
\begin{equation*}
\|(\Delta_{k} f)(\cdot, 0)\|^p_{L^p(\mathbb{R}^{n-1})}\lesssim 2^k
\int_{\mathbb{R}}\|M_{k, t}^{\ast}f(\cdot,
x_n)\|^p_{L^p(\mathbb{R}^{n-1})}dx_n,
\end{equation*}
Hence
\begin{equation}
\|(\Delta_{k} f)(\cdot, 0)\|_{L^p(\mathbb{R}^{n-1})}\lesssim
2^{k/p}\|M_{k}^{\ast}f\|_{L^p(\mathbb{R}^{n})}.  \label{max11}
\end{equation}
Write $\varphi'_{k}(\bar{x})$ as the BAPU functions in
$\mathbb{R}^{n-1}$. Then for fixed $k$, we have
\begin{align*}
(\mathscr{F}^{-1}_{\bar{\xi}}\varphi'_{k}\mathscr{F}_{\bar{x}})(\bar{x},
0)&=\sum_{l=k-1}^{\infty}(\mathscr{F}^{-1}_{\bar{\xi}}\varphi'_{k}\mathscr{F}_{\bar{x}}\mathscr{F}^{-1}\varphi_{l}\mathscr{F})(\bar{x},
0)\\
&=\sum_{l=k-1}^{\infty}(\mathscr{F}^{-1}_{\bar{\xi}}\varphi'_{k})\ast(\mathscr{F}^{-1}\varphi_{l}\mathscr{F})(\bar{x},
0)
\end{align*}

{\it Case 1.} $1\leq p\leq \infty$.  Using Young's inequality,
\eqref{max} and \eqref{max11}, we obtain
\begin{align*}
&\|\mathscr{F}^{-1}_{\bar{\xi}}\varphi'_{k}\mathscr{F}_{\bar{x}}f(\bar{x},
0)\|_{L^p(\mathbb{R}^{n-1})}\\
&\lesssim
\sum_{l=k-1}^{\infty}\|\mathscr{F}^{-1}_{\bar{\xi}}\varphi'_{k}\|_{L^1(\mathbb{R}^{n-1})}
\|\mathscr{F}^{-1}\varphi_{l}\mathscr{F}f\|_{L^p(\mathbb{R}^{n-1})}\\
&\lesssim \sum_{l=k-1}^{\infty} \|M_{l}^{\ast}f\|_{L^p(\mathbb{R}^{n})}\\
&\lesssim \sum_{l=k-1}^{\infty} 2^{l/p}\|\Delta_{l}
f\|_{L^p(\mathbb{R}^{n})}.
\end{align*}
Hence,
\begin{align*}
\|f(\bar{x}, 0)\|_{B^{s}_{p,q}(\mathbb{R}^{n-1})} \lesssim
\left(\sum_{k=0}^{\infty} 2^{skq}\left(\sum_{l=k-1}^{\infty}
2^{l/p}\|\Delta_{l} f\|_{L^p(\mathbb{R}^{n})}\right)^q
\right)^{1/q}.
\end{align*}
If $0<q \leq 1$, then
\begin{align*}
&\|f(\bar{x}, 0)\|_{B^{s}_{p,q}(\mathbb{R}^{n-1})} \lesssim
\left(\sum_{l=-1}^{\infty}\sum_{k=0}^{l+1}
 2^{skq}2^{lq/p} \|\Delta_{l}
f\|_{L^p(\mathbb{R}^{n})}^q \right)^{1/q}.
\end{align*}
If $s=0$, then
\begin{align*}
\|f(\bar{x}, 0)\|_{B^{s}_{p,q}(\mathbb{R}^{n-1})}& \lesssim
\left(\sum_{l=-1}^{\infty} l 2^{lq/p} \|\Delta_{l}
f\|_{L^p(\mathbb{R}^{n})}^q \right)^{1/q} \lesssim
\|f\|_{\tilde{B}_{ p,q}^{1/p}(\mathbb{R}^{n})}.
\end{align*}
In the case $s<0$,
\begin{align*}
\|f(\bar{x}, 0)\|_{B^{s}_{p,q}(\mathbb{R}^{n-1})} \lesssim
\left(\sum_{l=-1}^{\infty} 2^{lq/p} \|\Delta_{l}
f\|_{L^p(\mathbb{R}^{n})}^q\right)^{1/q}\lesssim \|f\|_{B_{
p,q}^{1/p}(\mathbb{R}^{n})}.
\end{align*}
 If $1\leq q \leq \infty$, using Minkowski's
inequality,
\begin{align*}
&\|f(\bar{x}, 0)\|_{B^{s}_{p,q}(\mathbb{R}^{n-1})} \lesssim
\sum_{l=-1}^{\infty}\left(\sum_{k=0}^{l+1} 2^{skq}2^{lq/p}
\|\Delta_{l}
f\|_{L^p(\mathbb{R}^{n})}^q  \right)^{1/q}\\
\end{align*}
then
\begin{align}
\|f(\bar{x}, 0)\|_{B^{s}_{p,q}(\mathbb{R}^{n-1})} \lesssim
\left\{\begin{array}{ll} \|f\|_{\tilde{B}_{
p,1}^{1/p}(\mathbb{R}^{n})}& s=0,\\
\|f\|_{B_{ p,1}^{1/p}(\mathbb{R}^{n})} &   s< 0.
\end{array}\right. \nonumber
\end{align}
{\it Case 2.} $0< p<1$.
\begin{align*}
&\|\mathscr{F}^{-1}_{\bar{\xi}}\varphi'_{k}\mathscr{F}_{\bar{x}}f(\bar{x},
0)\|^p_{L^p(\mathbb{R}^{n-1})}\\
&\lesssim
\sum_{l=k-1}^{\infty}2^{l(n-1)(1/p-1)p}\|\mathscr{F}^{-1}_{\bar{\xi}}\varphi'_{k}\|^p_{L^p(\mathbb{R}^{n-1})}
\|\mathscr{F}^{-1}\varphi_{l}\mathscr{F}f\|^p_{L^p(\mathbb{R}^{n-1})}\\
&\lesssim \sum_{l=k-1}^{\infty} 2^{l(n-1)(1-p)}2^{k(n-1)(p-1)}2^l\|M_{l}^{\ast}f\|^p_{L^p(\mathbb{R}^{n})}\\
&\lesssim \sum_{l=k-1}^{\infty}
2^{l(n-1)(1-p)}2^{k(n-1)(p-1)}2^l\|\Delta_{l}
f\|^p_{L^p(\mathbb{R}^{n})}.
\end{align*}
It follows that
\begin{align*}
\|f(\bar{x}, 0)\|_{B^{s}_{p,q}(\mathbb{R}^{n-1})} \lesssim
\left(\sum_{k=0}^{\infty} 2^{skq}\left(\sum_{l=k-1}^{\infty}
2^{l(n-1)(1-p)}2^{k(n-1)(p-1)}2^l\|\Delta_{l}
f\|^p_{L^p(\mathbb{R}^{n})}\right)^{q/p} \right)^{1/q}.
\end{align*}
If $q \leq p$, one has that
\begin{align*}
&\|f(\bar{x}, 0)\|^q_{B^{s}_{p,q}(\mathbb{R}^{n-1})} \\
&\lesssim \sum_{k=0}^{\infty} \sum_{l=k-1}^{\infty}
2^{l(n-1)(1/p-1)q}2^{k(s+(n-1)(1-1/p))q}2^{lq/p}\|\Delta_{l}
f\|^q_{L^p(\mathbb{R}^{n})} \\
&\lesssim \sum_{l=-1}^{\infty} \sum_{k=0}^{l+1}
2^{l(n-1)(1/p-1)q}2^{k(s+(n-1)(1-1/p))q}2^{lq/p}\|\Delta_{l}
f\|^q_{L^p(\mathbb{R}^{n})}.
\end{align*}
Recall that we write $s_{p} = (n-1)(1/(p\wedge 1) -1)$, then
\begin{align}
\|f(\bar{x}, 0)\|_{B^{s}_{p,q}(\mathbb{R}^{n-1})} \lesssim
\left\{\begin{array}{ll} \|f\|_{\tilde{B}_{ p,q}^{s_p+1/p
}(\mathbb{R}^{n})}& s=s_p,\\
\|f\|_{B_{ p,q}^{s_p+1/p}(\mathbb{R}^{n})}, & s<s_p.
\end{array}\right. \nonumber
\end{align}

 If $q \geq p$, using Minkowski's inequality, we have
\begin{align*}
&\|f(\bar{x}, 0)\|_{B^{s}_{p,q}(\mathbb{R}^{n-1})} \\
&\lesssim \left[\sum_{l=-1}^{\infty} \left( \sum_{k=0}^{\infty}
\left(\chi(k\leq
l)2^{l(n-1)(1-p)}2^{k(s+(n-1)(p-1))}2^{l}\|\Delta_{l}
f\|^p_{L^p(\mathbb{R}^{n})}\right)^{q/p}\right)^{p/q}\right]^{1/p} \\
&\lesssim \left[\sum_{l=-1}^{\infty} \sum_{k=0}^{l+1}
2^{l(n-1)(1/p-1)p}2^{k(s+(n-1)(1-1/p))p}2^{l}\|\Delta_{l}
f\|^q_{L^p(\mathbb{R}^{n})}\right]^{1/p}.
\end{align*}
Therefore,
\begin{align}
\|f(\bar{x}, 0)\|_{B^{s}_{p,q}(\mathbb{R}^{n-1})} \lesssim
\left\{\begin{array}{ll} \|f\|_{\tilde{B}_{p,p}^{s_p+1/p}(\mathbb{R}^{n})},& s=s_p,\\
\|f\|_{B_{ p,p}^{s_p+1/p}(\mathbb{R}^{n})}, &   s<s_p.
\end{array}\right. \nonumber
\end{align}

\rm In the case $1<p<\infty$, we define the maximum function
$M_{k,t}^{\ast} f$ as follows:
\begin{equation}
M_{k, t}^{\ast}f= \sup_{y\in \mathbb{Z}^n}\frac{|\Delta_{k,t}
f(x-y)|}{1+| 2^{k}y|^{n/r}}.\label{max5}
\end{equation}
Taking $y_1=... =y_{n-1}=0$, $y_n=x_n$ in \eqref{max5}, we have for
$ 2^{-k-1}\leq |x_n|\leq 2^{-k}$,
$$
|(\Delta_{k,t} f)(\bar{x}, 0)|\lesssim |M_{k, t}^{\ast}f(x)|, \quad
\bar{x}=(x_1, \cdots, x_{n-1})
$$
Hence
\begin{equation}
\|(\Delta_{k,t} f)(\cdot, 0)\|_{L^p(\mathbb{R}^{n-1})}\lesssim
\|M_{k, t}^{\ast}f(\cdot, x_n)\|_{L^p(\mathbb{R}^{n-1})},
\label{max7}
\end{equation}
Integrating \eqref{max7}, one has that
\begin{equation*}
\|(\Delta_{k,t} f)(\cdot, 0)\|^p_{L^p(\mathbb{R}^{n-1})}\lesssim 2^k
\int_{\mathbb{R}}\|M_{k, t}^{\ast}f(\cdot,
x_n)\|^p_{L^p(\mathbb{R}^{n-1})}dx_n,
\end{equation*}
Hence
\begin{equation}
\|(\Delta_{k,t} f)(\cdot, 0)\|_{L^p(\mathbb{R}^{n-1})}\lesssim
2^{k/p}\|M_{k, t}^{\ast}f\|_{L^p(\mathbb{R}^{n})}.  \label{max6}
\end{equation}
Let $\chi'_{k,t}(\bar{x})$ as the characteristic functions in
$\mathbb{R}^{n-1}$. Then for fixed $k$ and $t$, we have
\begin{align*}
(\mathscr{F}^{-1}_{\bar{\xi}}\chi'_{k,t}\mathscr{F}_{\bar{x}})(\bar{x},
0)&=\sum_{l=k}^{\infty}(\mathscr{F}^{-1}_{\bar{\xi}}\chi'_{k,t}\mathscr{F}_{\bar{x}}\mathscr{F}^{-1}\chi_{l,t}\mathscr{F})(\bar{x},
0)\\
&=\sum_{l=k}^{\infty}(\mathscr{F}^{-1}_{\bar{\xi}}\chi'_{k,t})\ast(\mathscr{F}^{-1}\chi_{l,t}\mathscr{F})(\bar{x},
0)
\end{align*}
Using Young's inequality, \eqref{max} and \eqref{max6}, we obtain
\begin{align*}
&\|\mathscr{F}^{-1}_{\bar{\xi}}\chi'_{k,t}\mathscr{F}_{\bar{x}}f(\bar{x},
0)\|_{L^p(\mathbb{R}^{n-1})}\\
&\lesssim
\sum_{l=k}^{\infty}\|\mathscr{F}^{-1}_{\bar{\xi}}\chi'_{k,t}\|_{L^1(\mathbb{R}^{n-1})}
\|\mathscr{F}^{-1}\chi_{l,t}\mathscr{F}f\|_{L^p(\mathbb{R}^{n-1})}\\
&\lesssim \sum_{l=k}^{\infty} \|M_{l,t}^{\ast}f\|_{L^p(\mathbb{R}^{n})}\\
&\lesssim \sum_{l=k}^{\infty} 2^{l/p}\|\Delta_{l,t}
f\|_{L^p(\mathbb{R}^{n})}.
\end{align*}
Hence,
\begin{align*}
\|f(\bar{x}, 0)\|_{B^{s}_{p,q}(\mathbb{R}^{n-1})} \lesssim
\left(\sum_{k=0}^{\infty} \sum_{t=1}^{T}
2^{skq}\left(\sum_{l=k}^{\infty} 2^{l/p}\|\Delta_{l,t}
f\|_{L^p(\mathbb{R}^{n})}\right)^q \right)^{1/q}.
\end{align*}
If $0<q \leq 1$, then
\begin{align*}
&\|f(\bar{x}, 0)\|_{B^{s}_{p,q}(\mathbb{R}^{n-1})} \\
&\lesssim \left(\sum_{l=0}^{\infty}\sum_{k=0}^{l-1}
\sum_{t=1}^{2^{n}} 2^{skq}2^{lq/p} \|\Delta_{l,t}
f\|_{L^p(\mathbb{R}^{n})}^q+\sum_{l=0}^{\infty}\sum_{k=l}
\sum_{t=2^{n}+1}^{T} 2^{skq}2^{lq/p} \|\Delta_{l,t}
f\|_{L^p(\mathbb{R}^{n})}^q  \right)^{1/q}.
\end{align*}
If $s=0$, then
\begin{align*}
\|f(\bar{x}, 0)\|_{B^{s}_{p,q}(\mathbb{R}^{n-1})}& \lesssim
\left(\sum_{l=0}^{\infty} \sum_{t=1}^{2^{n}}l 2^{lq/p}
\|\Delta_{l,t} f\|_{L^p(\mathbb{R}^{n})}^q+\sum_{l=0}^{\infty}
\sum_{t=2^{n}+1}^{T}
2^{lq/p} \|\Delta_{l,t} f\|_{L^p(\mathbb{R}^{n})}^q \right)^{1/q}\\
&\lesssim \|f\|_{\tilde{B}_{ p,q}^{1/p, 1/p}}.
\end{align*}
In the case $s<0$,
\begin{align*}
\|f(\bar{x}, 0)\|_{B^{s}_{p,q}(\mathbb{R}^{n-1})}& \lesssim
\left(\sum_{l=0}^{\infty} \sum_{t=1}^{2^{n}} 2^{lq/p} \|\Delta_{l,t}
f\|_{L^p(\mathbb{R}^{n})}^q+\sum_{l=0}^{\infty} \sum_{t=2^{n}+1}^{T}
2^{slq}2^{lq/p} \|\Delta_{l,t} f\|_{L^p(\mathbb{R}^{n})}^q \right)^{1/q}\\
&\lesssim \|f\|_{B_{ p,q}^{1/p,s+1/p}}.
\end{align*}
 If $1\leq q \leq \infty$, using Minkowski's
inequality,
\begin{align*}
&\|f(\bar{x}, 0)\|_{B^{s}_{p,q}(\mathbb{R}^{n-1})} \\
&\lesssim \sum_{l=0}^{\infty}\left(\sum_{k=0}^{l} \sum_{t=1}^{T}
2^{skq}2^{lq/p} \|\Delta_{l,t}
f\|_{L^p(\mathbb{R}^{n})}^q  \right)^{1/q}\\
&\lesssim \sum_{l=0}^{\infty}\left(\sum_{k=0}^{l-1}
\sum_{t=1}^{2^{n}} 2^{skq}2^{lq/p} \|\Delta_{l,t}
f\|_{L^p(\mathbb{R}^{n})}^q+\sum_{k=l} \sum_{t=2^{n}+1}^{T}
2^{skq}2^{lq/p} \|\Delta_{l,t} f\|_{L^p(\mathbb{R}^{n})}^q
\right)^{1/q}.
\end{align*}
then
\begin{align}
\|f(\bar{x}, 0)\|_{B^{s}_{p,q}(\mathbb{R}^{n-1})} \lesssim
\left\{\begin{array}{ll} \|f\|_{\tilde{B}_{
p,1}^{1/p, 1/p}(\mathbb{R}^{n})}& s=0,\\
\|f\|_{B_{ p,1}^{1/p, s+1/p}(\mathbb{R}^{n})} &   s< 0.
\end{array}\right. \nonumber
\end{align}

\medskip
$\hfill\Box$

\noindent \noindent{\bf Acknowledgment.} This work is partially
supported by the Marie Curie Excellence Project EUCETIFA,
MEXT-CT-2004-517154 of the European Commission. The second and third
named authors are supported in part by the National Science
Foundation of China, grant 10571004; and the 973 Project Foundation
of China, grant 2006CB805902.

\end{document}